\documentclass[11pt]{article}
\usepackage{amssymb} \usepackage{amsmath} \usepackage{amsthm}
\usepackage{latexsym} \usepackage{exscale}
\usepackage[final]{graphicx}
\usepackage[rflt]{floatflt}
\newtheorem{theorem}{Theorem}[section]
\newtheorem{lemma}[theorem]{Lemma}

\newtheorem{observation}{Observation}
 
\newenvironment{remark}{\stepcounter{remark}\noindent $\triangleleft$
 \textbf{Remark \thesection.\theremark: }}{$\triangleright$\newline}
 
\newcommand{\tit}{\textrm}

 \newcounter{remark}[section]

\def\eq#1{(\ref{#1})}

\def\({\left(\begin{array}{cccccc}}
\def\){\end{array}\right)}
\def\v{\vskip 1em}
\def\bes{\begin{eqnarray}}
\def\ees{\end{eqnarray}}
\def\forall{\hbox{for all~}}

\newcommand{\beq}{\begin{equation}}
\newcommand{\eeq}{\end{equation}}

\newcommand{\vp}{\varphi}
\newcommand{\ve}{\ensuremath{\varepsilon}} 
\newcommand{\si}{\ensuremath{\sigma}} 
\newcommand{\lam}{\ensuremath{\lambda}}

\newcommand{\eps}{\varepsilon}
 
\newcommand{\R}{\ensuremath{\mathbb{R}}}

\newcommand{\pf}{\begin{proof}}
\newcommand{\foorp}{\end{proof}}
\newcommand{\nquad}{\negthickspace\negthickspace
\negthickspace\negthickspace} 
\newcommand{\nqquad}{\nquad\nquad}

\numberwithin{equation}{section}

\begin{document}

\title{An Instability of the Godunov Scheme}
\author{Alberto Bressan\thanks{Department of Mathematics, Penn State
	University,
	UP, PA 16802, U.S.A.;
	E-mail: \texttt{bressan@math.psu.edu}}
	\and
	Helge Kristian Jenssen\thanks{Department of Mathematics, North
	Carolina State University,
	Raleigh, NC 27695, U.S.A.;
	E-mail: \texttt{hkjensse@math.ncsu.edu}. Research supported in part
	by the NSF under grant DMS-0206631}
	\and
	Paolo Baiti \thanks{Dipartimento di Matematica e Informatica,
	Universit\`a di Udine, Via delle Scienze 206,
	Udine 33100, Italy; E-mail: \texttt{baiti@dimi.uniud.it}}
	}
\date{\today}
\maketitle

\begin{abstract}
We construct a solution to a
$2\times 2$ strictly
hyperbolic system of conservation laws, showing that the
Godunov scheme \cite{Godunov59} can produce an arbitrarily large amount of oscillations.
This happens when the speed of a shock
is close to rational, inducing a resonance with the grid.
Differently from the Glimm scheme or the vanishing
viscosity method, for systems of conservation laws
our counterexample indicates that
no a priori BV bounds or $L^1$ stability estimates can in general
be valid for finite difference schemes.
\end{abstract}

\section{Introduction}\label{intro}
Consider a strictly hyperbolic $m\times m$ system of conservation laws
in one space dimension:
\beq\label{syst}
	u_t+f(u)_x=0.
\eeq
For initial data with small total variation,
the existence of a unique entropy weak solution
is well known \cite{Glimm65}, \cite{BressanLiuYang99}, \cite{BressanBook}.
A closely related question is the stability and convergence of various types
of approximate solutions.
For vanishing viscosity approximations
\beq\label{van_visc}
	u_t+f(u)_x=\ve\, u_{xx}\,,
\eeq
uniform BV bounds, stability and convergence as $\ve\to 0$
were recently established in \cite{BianchiniBressan04}.
Assuming that all the eigenvalues of the Jacobian matrix $Df(u)$
are strictly positive,
similar results are also proved in \cite{Bianchini03}
for solutions constructed by the semidiscrete (upwind)
Godunov scheme
\beq\label{semidiscrete}
	\frac{d}{dt}\, u_j(t)+ \frac{1}{ \Delta x}\Big[
	f\big(u_j(t)\big)-f\big(u_{j-1}(t)\big)\Big]~=~0\,,
	\qquad u_j(t)=u(t,~j\,\Delta x)\,.
\eeq
Recently, Bianchini has succeeded in proving the same type of
BV bounds for the Jin-Xin relaxation model \cite{Bianchini04}.

In the present paper we study the case of fully discrete
schemes, where the derivatives  with respect to both
time and space are replaced by
finite differences.  In sharp contrast with the previous situations,
we show that {\sl the total variation
of an approximate solution constructed by the Godunov scheme
can become arbitrarily large.}
This typically happens when the solution contains a shock
along a curve $x=\gamma(t)$, whose speed remains for a long time very
close (but not exactly equal) to a given rational, for example
\[\dot \gamma(t)\approx \frac{1}{ 2}\,\frac{\Delta x}{\Delta t}\,.\]
In this case, a resonance occurs and a substantial amount
of downstream oscillations is observed in the numerically computed solution.
For general $m\times m$ systems, our counterexample indicates that
a priori BV bounds
and uniform stability estimates cannot hold for solutions generated
by finite difference schemes.
It leaves open the possibility that these difference schemes
still converge to the unique entropy weak solution.  In any case,
a rigorous proof of this convergence cannot rely on
the same arguments as in \cite{Glimm65}, based on
uniform BV bounds and Helly's compactness theorem.

We recall that, for the $2\times 2$ system of isentropic gas dynamics,
the convergence of Lax-Friedrichs approximations was
proved in \cite{DingChenLuo89}, within the framework of compensated compactness.
Further results were obtained for {\sl straight line systems},
where all the Rankine-Hugoniot
curves are straight lines.
For these systems, uniform BV bounds, stability and convergence
of Godunov and Lax-Friedrichs approximations were established in \cite{LeVequeTemple85}, 
\cite{BressanJenssen00} and \cite{YangZhaoZhu03}.  The analysis relies on the fact that, due to the very
particular geometry, the interaction of waves of the same family
does not generate additional oscillations.

It is interesting to understand why the
arguments in \cite{BianchiniBressan04} or \cite{Bianchini03} break down, 
when applied to fully discrete
schemes.
A key ingredient in the analysis of vanishing
viscosity approximations is the local decomposition
of a viscous solution in terms of traveling waves.
To achieve a good control the oscillations produced by
interactions of waves of a same family, it is essential that
the center manifold of traveling profiles has a certain degree of smoothness.
This is precisely what fails in the case of fully discrete schemes.
As remarked by Serre \cite{Serre96}, for general hyperbolic systems
the discrete shock profiles cannot depend continuously on the speed $\sigma$,
in the BV norm.  In the related paper \cite{BaitiBressanJenssen04} we constructed 
an explicit example showing how this happens.
 
Our basic example is provided by a
$2\times 2$ system in triangular form
\begin{eqnarray}
	u_t+f(u)_x &=& 0,\label{u_eq}\\
	v_t+\lambda v_x+g(u)_x &=& 0 \label{v_eq}.
\end{eqnarray}
This system is strictly hyperbolic provided that
$f'(u)>\lambda$ for all $u$.
Choosing mesh sizes $\Delta t=\Delta x=1$,
the Godunov (upwind) scheme takes the form
\begin{eqnarray}
	u^{n+1}_j &=& u^n_j-\big[ f(u^n_j)-f(u^n_{j-1})\big],
\label{u_godunov}\\
	v^{n+1}_j &=& \lambda v^n_{j-1}+(1-\lambda) v^n_j -
	\big[g(u^n_j)-g(u^n_{j-1})\big].\label{v_godunov}
\end{eqnarray}
We shall  assume that
\beq\label{1.5}
	0<\lambda<f'(u)<1\qquad\qquad\forall ~u\in \R
\eeq
so that the usual linearized stability conditions are satisfied.
Thanks to the triangular form of the system \eq{u_eq}-\eq{v_eq}, 
the exact solution
of a Cauchy problem
can be computed explicitly. Indeed, one first solves the scalar equation for
$u$, say by the method of characteristics. Then the function
$g(u)_x$ is plugged as a source term into the second equation, 
which is linear
in $v$.

More specifically, we shall
consider a solution of \eq{u_eq}-\eq{v_eq} where the first component
$u$ contains a single shock, located along the curve
$ x=\gamma(t)$.   One can arrange things so that
the speed of the shock varies slowly in time,
remaining close to some rational number, say $\dot \gamma(t)\approx 1/2$.
The second component $v$ will then satisfy a linear
transport equation with a source located along a smooth curve,
and a priori bounds on its total variation can be easily given.

For the corresponding approximate solution generated by the
Godunov scheme \eq{u_godunov}-\eq{v_godunov}, however, things turn out to be quite different.
Because of the discretization,
in the equation for the second component
the source is not located along a smooth curve, but sampled at grid
points.  Since the speed of the shock is close to rational,
this source ``resonates'' with the grid, producing an
arbitrarily large amount of
downstream oscillations as time progresses.

The plan of the paper is as follows.
Section \ref{moving_source}
contains a preliminary analysis of
the heat equation with a moving source.
We show that if the point sources are
located at grid points and have an average speed
close to rational, then the solution
will contain downstream oscillations.
By studying the strength and the location of these tail oscillations
one gets a basic understanding of what happens for
solutions generated by discrete schemes.

In Sections \ref{special_solns} - \ref{total_variation} we carry out a detailed construction of a
Godunov approximate solution for the system \eq{u_eq}-\eq{v_eq},
showing that the total variation can become arbitrarily large.
Choosing the flux
\beq
	f(u)=\ln\big((1-\mu)+\mu e^u\big)\,,
\eeq
one can perform a nonlinear transformation introduced by Lax \cite{Lax54}, \cite{Lax57}
and explicitly compute the discrete Godunov solution $u^n_j$,
in terms of binomial coefficients.  In Section \ref{remarks} we summarize
the main features of our example, and discuss its significance
toward a rigorous theory of discrete schemes 
for hyperbolic conservation laws.
For readers' convenience, some results on the approximation 
of  the binomial distribution 
in terms of heat kernels are recalled in an appendix.

\section{The heat equation with a moving source}\label{moving_source}
As a motivation for the following analysis,
consider a solution of the finite difference
scheme
\beq\label{fin_diff}
	u^{n+1}_j=u^n_j-\big[ f(u^n_j)-f(u^n_{j-1})\big]
\eeq
in the form of a discrete shock with speed $\sigma$, say
\[u^n_j= \phi(j-\sigma n)\,.\]
For a scalar conservation law,
the existence of discrete traveling
profiles was proved in \cite{Jennings74}.
We assume that the flux $f$ is strictly convex and that the shock profile
satisfies
\[u^-~=~\lim_{s\to -\infty}\phi(s)~ >~\lim_{s\to +\infty}\phi(s)~=~u^+\,.\]
Inserting this solution as a source in the second equation \eq{v_godunov},
we obtain
\begin{eqnarray}
	&&\nqquad v^{n+1}_j = \lambda v^n_{j-1}+(1-\lambda) v^n_j -
	\big[g(u^n_j)-g(u^n_{j-1})\big] \nonumber\\
	&& = \lambda v^n_{j-1}+(1-\lambda) v^n_j -
	\Big[g\big(\phi(j-\sigma n)\big)-g\big(\phi(j-1-\sigma n)\big)
	\Big] \label{v_comp}.
\end{eqnarray}
We can choose a function $g$ such that
\begin{eqnarray*}
	g'(u) &=& 0\quad \mbox{if $\big|u-\frac{u^++u^-}{2} \big|\geq \epsilon$}\\
	g'(u) &>& 0\quad \mbox{if $\big|u-\frac{u^++u^-}{2}  \big|< \epsilon$,}
\end{eqnarray*}
for some $\epsilon>0$ small.
With this choice, the source terms on the right hand side of
\eq{v_comp} will vanish outside a thin strip centered around the shock.
Indeed, if the shock is located along the line
$x=\sigma t$, a fairly good approximation is
\beq
	g(u^n_j)-g(u^n_{j-1})~\approx~ \left\{
	\begin{array}{ll}
	g(u^+)-g(u^-)\qquad & \mbox{if $j-1=[\![\sigma n]\! ]\,$},\\
	0\qquad &  \mbox{otherwise}.
	\end{array}\right.
\eeq
Here $[\![s]\!]$ denotes the largest integer $\leq s$.
\begin{figure}
\centering
\includegraphics[width=9cm,height=7cm]
{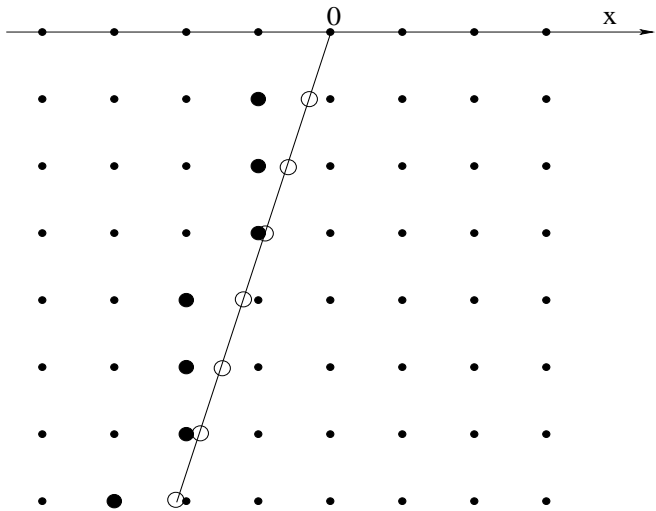}
\caption{Discrete sources at integer points.}
\end{figure} 
To understand the oscillations
in \eq{v_comp} determined
by the presence of these moving sources, we first study some model
problems related to the standard heat equation.
The approximation of a finite difference scheme by means of
a second order diffusion equation
is indeed a standard tool of analysis
(see \cite{LeVeque90}, p.117).   In this section
we review some calculations that highlight the mechanism 
responsible for generating variation in solutions to fully discrete 
schemes. For a detailed treatment we refer to the paper
\cite{BaitiBressanJenssen04}. 
Rigorous estimates will then be worked out in 
Sections \ref{special_solns} - \ref{total_variation}.

First, consider the easier case of the heat equation
with point sources located on a discrete set of points
$P_n=(n,\,\sigma n)$,
with $n$ integer
(the white circles in Figure 1)
\beq\label{delta_source}
	v_t-v_{xx}=\delta_{n,\sigma n}\,.
\eeq
We assume that $\sigma>0$ and  consider a solution of \eq{delta_source}
defined for $t\in ]-\infty,\,0[$.
Its values at time $t=0$ are now computed as
\beq\label{sum}
	v(0,y)=\Phi(y):=\sum_{n\geq 1} G(n,\, y+\sigma n)\,.
\eeq
Here 
\[G(t,x):= \frac{1}{2\sqrt{\pi t}}e^{-x^2/4t}\]
is the standard heat kernel.
To understand how the oscillations of $\Phi(y)$ decay as $y\to -\infty$,
we express the sum in \eq{sum} as an integral
\beq\label{Phi}
	\Phi(y):=\sum_{n\geq 1} G(n,\, y+\sigma n)=
	\int_0^\infty G(t,\,y+\sigma t)\, \big(1+h_1'(t)\big)\,dt
\eeq
where
\beq
	h_1(t):= [\![t ]\!]-t+1/2\,.
\eeq
By induction, we can find a sequence of periodic and uniformly bounded functions
$h_m$ such that
\[h_m(t)=h_m(t+1)\,,\qquad\int_0^1 h_m(t)\,dt=0\,,\qquad
\frac{d}{ dt}h_m(t)=h_{m-1}(t)\,.\]
Integrating by parts we obtain
\begin{eqnarray*}
	\Phi(y) &=& \int_0^\infty G(t,\,y+\sigma t)\,\Big(1+ \frac{d^m}{ dt^m}
	h_m(t)\Big)\,dt\\
	&=& \frac{1}{\sigma}+(-1)^m\int_0^\infty \frac{d^m}{ dt^m}
	G(t,\,y+\sigma t) \,h_m(t)\,dt\,.
\end{eqnarray*}
The identities
\[G(t,x)=t^{-1/2}G(1,~x/\sqrt t)\,,\qquad\qquad G_t=G_{xx}\,,\]
imply
\begin{eqnarray*}
	\frac{\partial^m}{\partial x^m}G(t,x) 
	&=& t^{-(m+1)/2}\cdot\frac{\partial^m G(1,~x/\sqrt t)}{\partial x^m},\\ \\
	\frac{\partial^m G(t,x)}{\partial t^m} 
	&=& t^{-(2m+1)/2}\cdot\frac{\partial^m}{\partial t^m}G(1,~x/\sqrt t)\,.
\end{eqnarray*}
In addition we observe that,
\[\sup_{\big|t+\textstyle\frac{y}{\sigma}\big|<|y|^{\epsilon+1/2}}
\left|\frac{d^m}{ dt^m}G(t,\,y+\sigma t) \right|~=~ O(1)\cdot e^{c_\epsilon y}\qquad
\qquad\hbox{as}~~y\to -\infty\,,\]
for some constant $c_\ve>0$.
Letting $y\to-\infty$, for every $m\geq 1$ the above estimates imply
\beq\label{estimate1}
	\left|\Phi(y)-\frac{1}{\sigma}\right| 
	\leq \int_0^\infty \left|\frac{d^m}{ dt^m}
	G(t,\,y+\sigma t) \right|\,dt~=~O(1)\cdot y^{-m/2}.
\eeq
Similarly,
\beq\label{estimate2}
	\big|\Phi'(y)\big|\leq \int_0^\infty \left|\frac{d^{m+1}}{ dt^{m+1}}
	G(t,\,y+\sigma t) \right|\,dt~=~O(1)\cdot y^{-(m+1)/2}.
\eeq
Taking $m=2$ in \eq{estimate2} one obtains the integrability of $\Phi'$,
hence a bound on the total variation of $\Phi$.

Next we outline the case when the point sources are located 
not at the points $P_n=(n,\sigma n)$, but at the points with 
integer coordinates
$Q_n:= \big(n,\,[\![\sigma n]\!]\big)$ (the black circles in Figure 1),
\[v_t-v_{xx}=\delta_{n,[\![\sigma n]\!]}\,.\]
Again we consider a solution defined for $t\in \,]-\infty,\,0]$
and a  direct computation yields
\[v(0,\,y-1)=\Psi(y):=\sum_{n\geq 1} G\big(n,\,y+[\![\sigma n]\!]\big).\]
Because of \eq{estimate1}, to determine the asymptotic behavior as $y\to-\infty$,
it suffices to estimate the difference
\[K(y)~:=~\Psi(y)-\Phi(y)~=~ -\sum_{n\geq 1} \Big[G(n,\,y+\sigma n\big)
-G\big(n,\,y+[\![\sigma n]\!]\big)\Big]\,.\]
It is here that, if the speed $\sigma$ is close to a rational,
a resonance is observed.   To see a simple case, let
$\sigma=1+\ve$, with $\ve>0$ small.  Then
we can approximate
\begin{eqnarray}\label{K}
	K(y)&\approx~& -\sum_{n\geq 1} G_x(n,\,y+\sigma n\big)
	\big(\sigma n-[\![\sigma n]\!]\big)\nonumber \\
	&\approx ~& -\int_0^\infty G_x(t,\,y+\sigma t)\, \big(
	\ve t -[\![\ve t]\!]\big)\,dt\,.
\end{eqnarray}
\begin{figure}
\centering
\includegraphics[width=9cm,height=6cm]
{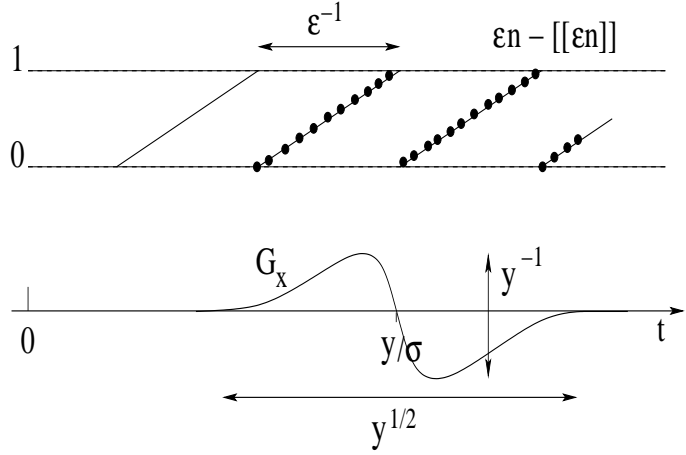}
\caption{Interaction of $G_x$ and fractional parts.}
\end{figure} 
The functions appearing in the above integration are shown in Figure 2.
We recall that
\[\int_0^\infty G_x(t,\,y+\sigma t)\,dt~=~
-\int_0^\infty \frac{y+\sigma t}{ 4t\sqrt{\pi t}}\,
\exp\left\{ -\frac{(y+\sigma t)^2}{ 4t}\right\}\,dt~=~0\]
for every $y<0$. Set $y_\ve:= -\ve^{-2}$.
When $y$ ranges within the interval
\[I_\ve:= [y_\ve,~y_\ve/2]~=~[-\ve^{-2},~-\ve^{-2}/2]\,,\]
the integral in \eq{K} can be of the same order of magnitude as
\[\int_0^\infty \big|G_x(t,\,y_\ve+\sigma t)\big|\,dt~\geq~
c_0\, y_\ve^{-1/2}~=~c_0\,\ve\,.\]
Moreover, each time that
$y$ increases by an amount $\Delta y=\ve^{-1}$, the phase of the
fractional part $[\![\ve y]\!]-\ve y$ goes through a full cycle, hence the map
\[y\mapsto \int_0^\infty G_x(t,\,y+\sigma t)\, \big(\ve t -[\![\ve t]\!]\big)\,dt\]
oscillates by an amount $\geq c_1\ve$.  In all, we have approximately
$1/2\ve$ cycles within the interval $I_\ve$. Hence the total variation of the
discrete profile $\Psi=\Psi^{(1+\ve)}$ can be estimated as
\beq\label{TV}
	\mbox{Tot. Var. }\big\{
	\Psi^{(1+\ve)}~;~~I_\ve \big\}\geq c_2 
\eeq
for some constant $c_2>0$ independent of $\ve$.
By \eq{TV} it is clear that, as $\ve\to 0+$, the
functions $\Psi^{(1+\ve)}$
do not form a Cauchy sequence and cannot converge in the space BV.
The details of the preceding argument are given (for the Lax-Friedrichs scheme) 
in \cite{BaitiBressanJenssen04}.
\begin{figure}
\centering
\includegraphics[width=9cm,height=7cm]
{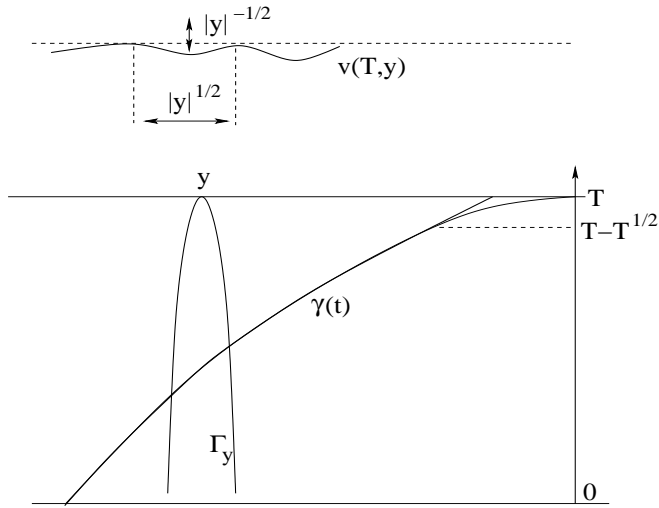}
\caption{Oscillations produced by sources along the 
curve $\gamma(t)$.}
\end{figure} 
Finally we consider the case where the source travels with a 
carefully chosen variable speed. In the previous 
computation the source travelled with a constant
speed $\sigma=1+\ve$, and a significant amount of
oscillation was observed at a distance
$|y_\ve|=\ve^{-2}$ downstream from the shock.
To construct a solution whose oscillation becomes
arbitrarily large, the
idea is to choose a large time interval $[0,T]$ and
generate a source with non-constant speed,
say located along a curve $x=\gamma(t)$ with
\beq\label{gamma_speed}
	\gamma(t)=(t-T)-2\sqrt{T-t}\,,\qquad
	\dot \gamma(t)= 1+\frac{1}{ \sqrt{T-t}}
\eeq
for $0\leq t\leq T-T^{1/2}$, (see Figure~3).
Consider the profile at the terminal time $T$
of a solution of
\beq
	v_t+v_{xx}=\delta_{n,\,[\![\gamma(n)]\!]}\qquad\qquad v(0,x)=0,
\eeq
defined for $t\in [0,T]$.
For a given point $y$ the value $v(T,y)$ is 
essentially determined by the sources
located inside the parabolic region
\[\Gamma_y:=\Big\{ (t,x)~:~~x\in \big[ y-(T-t)^{\epsilon+1/2},~
	y+(T-t)^{\epsilon+1/2}\big]\Big\}.\]
By construction, if $-T/2<y<- T^{1/2}$, this region contains a
portion of the curve $\gamma$
traveling with speed
\[\dot \gamma(t)=1+\frac{1}{ \sqrt{T-t}}\approx 1+\frac{1}{ \sqrt{|y|}}\,.\]
This is precisely the resonant speed that can produce downstream
oscillations at a distance $\approx y$ from the shock.
As in \eq{TV}, one has
\beq\label{2.14}\mbox{Tot. Var. }\big\{ v(T,\cdot)~
;~~[y,\, y/2] \big\}\geq c_2 \qquad\qquad -T/2<y<- T^{1/2}\,.
\eeq
Now consider the points
\[ y_j\doteq - 2^j\,T^{1/2}\qquad\qquad j=0,1,\ldots,N-1\,,\]
where $N$ is the largest integer $\leq \frac{1}{2}\log_2 T $.
Notice that this choice of $N$ guarantees that
$-T/2<y_N<\ldots<y_1<y_0=-T^{1/2}$.
Observe that the choice of the variable speed $\dot \gamma$
at \eq{gamma_speed} is precisely what is needed
in order to produce a uniformly positive
variation on each of the intervals 
$[y_j,\, y_{j-1}]$. From \eq{2.14} it thus follows
\beq\label{2.15}\mbox{Tot. Var. }\big\{ v(T,\cdot)\big\}\geq 
\sum_{j=1}^N \mbox{Tot. Var. }\big\{ v(T,\cdot)~
;~~[y_j,\, y_j/2] \big\}\geq c_2\cdot\left(\frac{1}{2}\log_2 T - 2
\right)\,.
\eeq
Choosing $T$ large, we expect to
find numerical
solutions generated by the Godunov scheme
whose total variation grows by an arbitrarily large
amount.  In the following sections we show that this
is indeed the case. By providing rigorous estimates on all
the approximations performed in the above formal analysis,
we will prove that the discrete Godunov solution
satisfies the same type of estimate as in \eq{2.15}.

\section{Special solutions}\label{special_solns}
\subsection{Inviscid solution}\label{inviscid_soln}
It will be convenient to work with the particular scalar conservation law
\beq\label{inviscid}
	u_t+\big[\ln(\nu  +\mu  e^u )\big]_x=0\,,
\eeq
where $\mu \in\, ]\,1/2,1\, [$ and $\nu:=1-\mu$. This equation will serve as the first equation in the 
system \eq{u_eq}-\eq{v_eq}. Notice that in this case the flux 
$f(u)=\ln(\nu  +\mu e^u)$ is increasing and convex. Indeed,
\[f'(u)=\frac{\mu e^u}{ \nu  +\mu e^u}\,, \qquad f''(u)=\frac{\mu\nu  e^u}{ (\nu  +\mu e^u)^2}\,,\]
so that
\[f'(u)\in\, ]\,0,1\, [\,,\qquad\qquad f''(u)>0.\]
Note that 
\[f(0)=0,\qquad\qquad f'(0)=\mu>1/2.\]
Our first goal is to construct an exact solution
of the inviscid equation \eq{inviscid} containing exactly one shock, located
along a curve $\gamma(t)$ which travels with a speed close to a rational 
number $\lam_0>\lam$, where $\lam$ is as in \eq{v_eq}. 
The role of  this curve will be that its spatial translates are 
close to the level curves of a certain (everywhere defined) solution of the 
upwind scheme \eq{u_godunov}. This particular solution of the scheme 
will then be a good approximation of the exact solution of \eq{inviscid} which
we construct here.

We now fix the 
speed $\lam$ in \eq{v_eq} to be $1/2$. (This choice is not essential but it 
simplifies the calculations below which involve approximating binomial 
distributions with the heat kernel.) The shock solution of \eq{u_eq} will 
connect a variable state $u_-(t)$ to a fixed right state $u_+<u_-(t)$, and 
we want this $u$-shock to travel faster than the advection speed $\lambda=1/2$ 
of the second equation. We begin by fixing the right state $u_+:=0$, and 
our solution will be constant equal to this value to the right of $\gamma$, i.e.
\[u(t,x)\equiv 0\qquad\qquad \mbox{for $x>\gamma(t)$}.\]

Next, in analogy with the analysis in Section \ref{moving_source}, we insist
that the shock curve $\gamma(t)$ should have a speed as in \eq{gamma_speed}. 
For notational convenience we will make the explicit construction for 
negative times. We thus fix a large, positive (integer) time $T$, 
and define the curve $\gamma(t)$ by setting
\beq\label{gamma}
	\gamma(t)=\lam_0 t-2\sqrt{-t},
	\qquad  \mbox{for $t\in\, [\,-2T,-\sqrt{T}\, [\, $}. 
\eeq
We require $\gamma$ to be smooth on $[-2T,0]$ and such that 
\[\dot \gamma (t) \equiv \dot\gamma(-\sqrt{T}) = \lambda_0+\frac{1}{T^{1/4}}
\qquad  \mbox{for $t\in\, [\,-\sqrt{T},0\, ]\, , $}\]
see Figure 3. Here $\lam_0\in\, ]\,\mu,1\, [\, \cap\, \mathbb Q$ denotes a fixed, rational  speed. 
At every time $t\in \, [\, -2T,0\, ]$ the state $u_-(t)$ immediately to the left of the shock 
must then satisfy the Rankine-Hugoniot equation
\beq\label{RH1}
	\frac{\ln(\nu  +\mu e^{u_-(t)})}{u_-(t)}=\dot \gamma(t).
\eeq
To establish the existence and properties of the function
$t\mapsto u_-(t)$, we let $u_0$ denote the left state corresponding to 
constant speed $\lam_0$, i.e. 
\[\frac{\ln(\nu  +\mu e^{u_0})}{u_0}=\lam_0 \,,\] 
and then consider the equation
\beq\label{RH2}
	\frac{\ln(\nu  +\mu e^u)}{u}=\sigma \,.
\eeq
A straightforward argument shows that \eq{RH2} has exactly one solution
$u=\vp(\si)\in\, ]\,0,\infty\,[$ for each $\si\in\, ]\, \mu,1\, [$, and that $\vp$ depends smoothly on $\si$.
Therefore, for $t$ negative and sufficiently large, from
\eq{RH1} and \eq{gamma} we recover a smooth function $t\mapsto u_-(t)$, where
\begin{eqnarray}
	u_-(t) &=& \vp\left( \lam_0+\frac{1}{\sqrt{-t}}\right)\nonumber \\
	&=& \kappa_0+\frac{\kappa_1}{ (-t)^{1/2}}+\frac {\kappa_2}{(-t)}
	+\frac{O(1)}{ (-t)^{3/2}}\,, \label{taylor1}
\end{eqnarray}
for $t<0$. We have here assumed the Taylor expansion
\[\vp(\lam_0+\delta)=\kappa_0 + \kappa_1\delta
	+ \kappa_2\delta^2+O(1)\delta^3.\]

\subsection{Discrete Cole-Hopf transformation}\label{discr_appr}
As observed by Lax \cite{Lax54}, \cite{Lax57},
for the particular flux function $f$ in \eq{inviscid}, one can 
perform a ``discrete Cole-Hopf" transformation which linearizes 
the scheme \eq{u_godunov}. This enables us to write down  
explicit  solutions of the Godunov scheme in the form of 
nonlinear superpositions of discrete traveling waves. More precisely, 
consider the Godunov scheme in this case, i.e.
\beq\label{spec_godunov}
	u^{n+1}_j=u^n_j-\Big[ \ln\big(\nu  +\mu e^{u^n_j}\big)-
	\ln\big(\nu  +\mu e^{u^n_{j-1}}\big)\Big],
\eeq
where we recall that $\mu>1/2$ and that $\nu=1-\mu$. We then have the following
result due to Lax.
\begin{lemma}\label{Cole-Hopf}
If $z^n_j>0$ is a solution of
\beq\label{linearized_scheme}
	z^{n+1}_j=\mu z^n_{j-1}+\nu  z^n_j\,,
\eeq
then a solution of {\rm\eq{spec_godunov}} is provided by
\[u^n_j=\ln\left( \frac{z^n_{j-1}}{ z^n_j}\right)\,.\]
\end{lemma}
\begin{proof} A direct computation yields
\begin{eqnarray*}
	u^{n+1}_j &=& \ln\left(\frac{ z^{n+1}_{j-1}}{ z^{n+1}_j}\right)
	=\ln\left( \frac{\mu z^n_{j-2}+\nu  z^n_{j-1}}{ \mu z^n_{j-1}+\nu  z^n_j}\right)\\ \\
	&=&  \ln \left(\frac{ z^n_{j-1}}{ z^n_j}\right)+\ln \left(\frac{ \nu   +\mu \frac{z^n_{j-2}}{ z^n_{j-1}}}
	{\nu   +\mu \frac{z^n_{j-1}}{ z^n_j}}\right) \\ \\
	&=&u^n_j+ 
	\ln\left( \nu  +\mu e^{u^n_{j-1}}\right)-
	\ln\left(  \nu  +\mu e^{u^n_j}\right)\,.
\end{eqnarray*}
\end{proof}

\noindent In this section we will give discrete traveling waves 
for the the linearized scheme \eq{linearized_scheme} and use these
to construct discrete traveling waves for \eq{spec_godunov}. In the next section 
we will then use these as building blocks to write down an explicit solution of 
\eq{spec_godunov} which is a  good approximation to the exact 
inviscid solution described in Section \ref{inviscid_soln}.
Letting both variables in \eq{linearized_scheme} range over $\R$ we get
\beq\label{cont_linearized}
	z(t+1,x)=\mu z(t,x-1)+ \nu   z(t,x).
\eeq
Special solutions of \eq{cont_linearized} are easy to find.
For every $b>0$ define $\sigma(b)$ as the Rankine-Hugoniot
speed of the shock connecting the left state $b$ to the right state $u_+=0$, 
i.e.
\beq\label{speed}
	\sigma(b):= \frac{\ln (\nu  +\mu e^b)}{b}\,.
\eeq
Notice that $\sigma=\sigma(b)$ then satisfies the equation
\beq\label{relation}
	\nu   + \mu e^b=e^{b\si(b)},
\eeq
and it follows that for any $x_1\in\R$ the function
\[z(t,x)=e^{-b[x-x_1-\sigma(b) t]}\]
is a solution of the linearized equation \eq{cont_linearized}.  
More generally, by linearity, any finite linear combination
\[z(t,x)= \sum_{j=0}^N e^{-b_j[x-x_j-\sigma(b_j)t]},\]
or any integral combination
\beq\label{int_comb}
	z(t,x)=\int e^{-\zeta\big(x-x(\zeta)-\sigma(\zeta)t\big)}\,d\zeta,
\eeq
are also solutions of \eq{cont_linearized}.  In turn, the formula
\beq\label{lax-cole-hopf}
	u(t,x)=\ln \left( \frac{z(t,\,x-1)}{ z(t,x)}\right)
\eeq
allows us to recover an explicit solution of the nonlinear 
difference equation
\beq\label{u_cont}
	u(t+1,\,x)=u(t,x)-\Big[ \ln\big(\nu  +\mu e^{u(t,x)}\big)-
	\ln\big(\nu  +\mu e^{u(t,x-1)}\big)\Big]\,.
\eeq
As an example let $b_1<b_2$, $x_1,x_2\in\R$ and consider the
solution 
\[z(t,x):= e^{-b_1[x-x_1-\sigma(b_1) t]}+
e^{-b_2[x-x_2-\sigma(b_2) t]}\]
of \eq{cont_linearized}.
A simple calculation shows that the corresponding 
function $u$ in \eq{lax-cole-hopf}
represents a traveling wave solution of \eq{u_cont} with speed 
\[\sigma^*=\frac{\ln(1+e^{b_1})-\ln(1+e^{b_2})}{ b_1-b_2}=\frac{b_1\sigma(b_1)-
	b_2\sigma(b_2)}{ b_1-b_2},\]
and connecting $u(-\infty)=b_2$ to $u(+\infty)=b_1$.
To see where the center of the wave is located, we observe that
at the point
\[\bar x=\frac{x_1b_1-x_2b_2}{ b_1-b_2}= x_2+\frac{b_1}{ b_1-b_2}(x_1-x_2)\]
one has
\[u(0, \bar x)=\ln\left( \frac{e^{b_1} e^{-b_1(\bar x-x_1)}+
e^{b_2} e^{-b_2(\bar x-x_2)}}{
e^{-b_1(\bar x-x_1)}+
e^{-b_2(\bar x-x_2)}}\right)=\ln\left(\frac{e^{b_1}+e^{b_2}}{ 2}\right)\approx
\frac{b_1+b_2}{ 2}\,.\]
We can thus think of the point $x^*(t):=
\bar x+\sigma^*t$ as the center of the discrete traveling profile at time $t$.

\section{Construction of a special discrete approximation}
\label{special_discr_appr}
We proceed to use the information from Section 
\ref{discr_appr} to construct an exact solution of the upwind 
scheme \eq{spec_godunov} that is  a good approximation 
of the particular solution of \eq{inviscid} with a single shock 
along the curve $\gamma$.  This exact solution 
will be obtained as a solution $u^d (x,t)$ of the difference 
equation \eq{u_cont} as in \eq{lax-cole-hopf}, where $z^d(t,x)$ 
is a special solution of the linearized 
difference equation \eq{cont_linearized}. (Superscript $d$ 
is used (only) here to emphasize that the solutions satisfy 
the everywhere defined difference equations \eq{u_cont}
and \eq{cont_linearized}.)
We point out that the relevant property of the  ``good" 
approximation $u^d$, which will be needed in the final 
computations, is that its level curves
should travel with speeds sufficiently close to the level curves 
of the exact inviscid solution of \eq{inviscid} constructed in 
Section \ref{inviscid_soln}. 

To construct such a solution of \eq{u_cont} we consider the curve
$\gamma(t)$ given by \eq{gamma} for large negative times. As above we set
$u_-(t):=\vp(\dot\gamma(t))$, such that $\sigma(u_-(\tau))=\dot\gamma(\tau)$, 
and define
\beq\label{z_1}
	z(t,x) := 1+ \int_{-2T}^{\,0}\!\!\! \exp\Big\{-u_-(\tau)
	\,\big[x-\gamma(\tau)-\dot\gamma(\tau)(t-\tau)\big]\Big\}
	\cdot\psi(\tau)\,d\tau\,,
\eeq
for $t\in\mathbb R$. 
By \eq{int_comb} this defines a smooth solution of the difference equation 
\eq{cont_linearized}. The corresponding solution of \eq{u_cont} will be 
denoted by $u(t,x)$ (see \eq{U_def} below). For reasons that will be clear from the 
following computations we choose 
\beq\label{psi}
	\psi(\tau):=(-\tau)^{-3/4},\qquad\qquad \tau< 0.
\eeq
To avoid several minus signs we let $s=-t$, $\eta=-x$, and set
\begin{eqnarray}
	&&Z(s,\eta) := z(t,x) \label{Z}\\
	&&\qquad = 1+\int_{0}^{2T}
	\exp\Big\{ u_-(-\tau)\big[ \eta+\gamma(-\tau)-\dot\gamma(-\tau)(s-\tau)\big] \Big\}
	\psi(-\tau)\, d\tau, \nonumber
\end{eqnarray}
where we have made the change of variables $\tau\mapsto -\tau$. We define 
\begin{eqnarray}
	\Gamma(\tau) &:=& -\gamma(-\tau)=\lam_0\tau +2\sqrt{\tau},\label{Gamma}\\
	W(\tau) &:=& u_-(-\tau)=\varphi\left(\lam_0+\frac{1}{\sqrt{\tau}}\right),\label{W}\\
	g(\tau) &:=& \psi(-\tau)=\tau^{-3/4}, \label{g}
\end{eqnarray}
for $\tau>0$.
Letting $y:=\eta-\Gamma(s)$ and making the change of variables $\xi=\tau/s$, we obtain
\begin{eqnarray}
	&&\nqquad\nqquad Z(s,\Gamma(s)+y) \nonumber \\ 
	&&\nquad\nqquad=  1+\int_{0}^{2T} 
	\exp\left\{ W(\tau)\left[y+\Gamma(s)-\Gamma(\tau) -\dot\Gamma(\tau)(s-\tau) \right] \right\}
	g(\tau)\, d\tau\nonumber \\ \nonumber \\
	&&\nquad\nqquad= 1+\int_{0}^{2T}
	\exp\left\{ W(\tau)\left[ y +\sqrt{s}\left(2-\sqrt{\frac{\tau}{s}}-\sqrt{\frac{s}{\tau}}\right) \right] \right\}
	g(\tau)\, d\tau\nonumber \\ \nonumber \\
	&&\nquad\nqquad= 1+s^{1/4}\int_{0}^{2T/s}
	\exp\left\{ W(s\xi)\Big[ y +\sqrt{s}\Big(2-\sqrt{\xi}-\frac{1}{\sqrt{\xi}}\Big) \Big] \right\}
	g(\xi)d\xi\, .\label{int}
\end{eqnarray}
We will consider times $t$ in $[\,-T,0\,]$, i.e. $s\in\, [\,0,T\,]$, such that the upper limit of integration in 
\eq{int} is bounded away from $\xi=1$. The factor $s^{1/4}$ will be important 
in the following computations and this is what  dictates the particular choice in \eq{psi}. 
We proceed to analyze in detail the level curves of $z(t,x)$ and $u(t,x)$ by considering 
the level curves of $Z(s,\eta)$. 
For $s\gg 1$, the main contribution to the integral in \eq{int}
is provided within the region where $\xi\approx 1$. A careful analysis
is required in order to obtain the exact behavior.

\subsection{Level curves of $u(t,x)$}\label{u_level}
As a first step we show that the time derivatives of the level 
curves of $Z(s,\eta)$ behave (almost) like  $\dot\Gamma(s)$.
Before starting to estimate the various terms that are needed, we make 
a simplifying observation.
For $\eps>0$, $\eps\ll 1$, we let $I(s)$ denote the interval
\[I(s)=[1-s^{-\alpha},1+s^{-\alpha}], \qquad \qquad \alpha =1/4-\eps.\]
\begin{observation}\label{exp}
	In what follows all $\xi$-integrals over $[\, 0,2T/s\,]\setminus I(s)$, 
	as well as all boundary terms (obtained by 
	differentiation of such integrals), are exponentially small.
	Namely, they are all of order $O(e^{-Cs^\eps})$ for some $C,\, \eps>0$.
\end{observation}
\noindent Next we define the functions
\begin{eqnarray}
	F(\xi,s,y) &:=&  W(s\xi)\left[ y +\sqrt{s}\left(2-\sqrt{\xi}-\frac{1}{\sqrt{\xi}}\right) \right] ,\label{F} \\ \nonumber \\
	Q(s,y) &:=& \left(\kappa_0 + \frac{\kappa_1}{\sqrt{s}}\right)y, \label{Q}\\ \nonumber \\
	H(s,\xi) & := & -\frac{\kappa_0\sqrt{s}}{4}(\xi-1)^2.\label{H}
\end{eqnarray}
Thus 
\[Z(s,\Gamma(s)+y)=1+s^{1/4}\int_{0}^{2T/s} \exp\big[F(\xi,s,y)\big]g(\xi)\, d\xi.\]
\begin{lemma}\label{dZ/ds}
	For $y=O(1)$ and $s\gg 1$ we have 
	\beq\label{Z'}
	\frac{d}{ds}Z(s,\Gamma(s)+y)=\frac{O(1)}{s^{3/2-\eps}},
	\eeq
\end{lemma}
\begin{proof}
We start by stating a detailed expansion of the function $F-Q-H$
where $F$, $Q$ and $H$ were defined by \eq{F}, \eq{Q} and \eq{H}:
\begin{eqnarray*}
	F(\xi,s,y) -Q(s,y)-H(s,\xi) & = & A+B(\xi-1)+C(\xi-1)^2+ D(\xi-1)^3\\ 
	&&+E(\xi-1)^4 +O(1)(\xi-1)^5,
\end{eqnarray*}
	where the coefficients $A$-$E$ are given as
\begin{eqnarray*}
	A(y,s)&=&y\left( \frac{\kappa_2}{s} + \frac{O(1)}{s^{3/2}}\right),\\ 
	B(y,s)&=& -y\left( \frac{\kappa_1}{2\sqrt{s}} + \frac{\kappa_2}{s} +\frac{O(1)}{s^{3/2}}\right),\\ 
	C(y,s)&=& -\frac{\kappa_1}{4} + \frac{3\kappa_1y-2\kappa_2}{8\sqrt{s}} 
		+ \frac{\kappa_2y+O(1)}{s} + \frac{O(1) y}{s^{3/2}},\\
	D(y,s) &=&  \frac{\kappa_0\sqrt{s}}{4} + \frac{3\kappa_1}{8} 
		+ \frac{8\kappa_2-5y\kappa_1}{16\sqrt{s}}+ \frac{O(1)-\kappa_2 y}{s} 
		+ \frac{O(1)y}{s^{3/2}},\\
	E(y,s) &=&O(1)\sqrt{s}+\frac{O(1)y}{\sqrt{s}} .
\end{eqnarray*}
The proof is a direct calculation which is omitted.
The estimate \eq{Z'} is now obtained by exploiting the fact
that, when differentiating with respect to time, we can subtract off the leading 
order contribution, i.e. the term
\[\frac{d}{ds}\left\{ s^{1/4}\int_{I(s)} e^{H(s,\xi)}\, d\xi \right\},\]
since this is exponentially small.
We have
\[Z(s,\Gamma(s)+y)=1+e^{Q(s,y)} s^{1/4}\int_{I(s)}\!\!\! \exp[F(\xi,s,y)-Q(s,y)]g(\xi)\, d\xi.\]
such that, up to exponentially small terms, 
\begin{eqnarray}
	\frac{d}{ds} Z(s,\Gamma(s)+y)&=& Q_se^{Q}\cdot s^{1/4}\int_{I(s)} \exp[F-Q]g(\xi)\, d\xi   \nonumber\\ \nonumber\\
	&&+e^{Q}\frac{\partial}{\partial s}\left\{ s^{1/4}\int_{I(s)}
	\exp[F-Q]g(\xi)-e^{H}\, d\xi \right\}\nonumber\\ \nonumber\\
	&&\nqquad\nqquad \nquad
	=\frac{O(1)}{s^{3/2-\eps}} + \frac{O(1)}{s} \cdot s^{1/4}\int_{I(s)}\big(\exp[F-Q-H]g(\xi)-1\big)e^{H}\, d\xi  
	\nonumber\\ \nonumber\\
	&&\nqquad\nqquad + O(1)s^{1/4}\int_{I(s)} \Big\{\big(\exp[F-Q-H]g(\xi)-1\big)e^{H}\Big\}_s\, d\xi, \label{first}
\end{eqnarray}
Using the expansion above and $g(\xi)=1-3(\xi-1)/4 +O(1)(\xi-1)^2$, we get that (for $y=O(1)$ and 
$s\gg 1$)
\beq\label{F-Q-H}
	\exp[F-Q-H]g(\xi)-1= \frac{O(1)}{s^{1/2-\eps}} + \{\mbox{odd function of $(\xi-1)$}\}.
\eeq
For the last term in \eq{first} we apply the expansion above again to get that 
\begin{eqnarray*}
	&&\nquad\Big\{\big(\exp[F-Q-H]g(\xi)-1\big)e^{H}\Big\}_s \\ \\
	&&\nquad=\Big[\big(\exp[F-Q-H]g(\xi)-1\big)H_s
	+\exp[F-Q-H]g(\xi)(F-Q-H)_s\Big] e^H \\ \\
	&&\nquad = \left[\{\mbox{odd function of $(\xi-1)$}\}+\frac{O(1)}{s^{3/2-\eps}}\right] e^{H}.
\end{eqnarray*}
Here we have used the expansion of $g(\xi)$, \eq{F-Q-H}, 
that 
\[H_s=-\frac{\kappa_0(\xi-1)^2}{8\sqrt{s}}=\frac{O(1)}{s^{1-\eps}}\qquad\mbox{on $I(s)$,} \]
and finally that
\[(F-Q-H)_s=\frac{\kappa_0}{2\sqrt{s}}(\xi-1)^3+\frac{O(1)}{s^{3/2-\eps}},\]
which follows from the expansion above. As
\[s^{1/4}\int_{I(s)} e^{H(s,\xi)}\, d\xi=O(1),\]
we therefore have that the two last terms in \eq{first} are both $O(1)/s^{3/2-\eps}$. 
The result  follows.
\end{proof}

\noindent We now define 
\beq\label{U_def}
	U(s,\eta):=u(-s,-\eta):=\ln\left(\frac{z(-s,-\eta-1)}{z(-s,-\eta)}\right)=\ln\left(\frac{Z(s,\eta+1)}{Z(s,\eta)}\right).
\eeq
A similar calculation as above shows that
\[Z(s,\Gamma(s)+y)=1+ c_0e^{Q(s,y)} + \frac{O(1)}{s^{1/2-\eps}},\]
where 
\[c_0=t^{1/4}\int_{-\infty}^\infty\exp[H(t,\xi)]\, d\xi=2\sqrt{\frac{\pi}{\kappa_0}}\, .\]
It follows that 
\[U(s,\Gamma(s)) = \ln\left(\frac{1+c_0e^{\kappa_0}}{1+c_0}\right)+ 
\frac{O(1)}{s^{1/2-\eps}}=:a_0 + \frac{O(1)}{s^{1/2-\eps}}.\]
For $a$ sufficiently close to $a_0$, and well within $]\, U(-\infty), U(+\infty)\,[$, we 
denote the $a$-level curve of $U(s,\eta)$ by $X(s;a)$, i.e.
\beq\label{U_level}
	U(s,X(s;a)) \equiv a,
\eeq
such that  
\beq\label{X_s}
	U_s(s,X(s;a))=-X_s(s;a) U_\eta(s,X(s;a)).
\eeq
A computation shows that 
\[Z_\eta(s,\Gamma(s)+y)=\kappa_0c_0e^{Q(s,y)} + \frac{O(1)}{s^{1/2-\eps}},\]
whence 
\[U_\eta(s,\Gamma(s)+y)=\frac{\kappa_0c_0e^{\kappa_0 y}(e^{\kappa_0}-1)}
{(1+\kappa_0c_0e^{\kappa_0y})(1+c_0e^{\kappa_0y})}+ \frac{O(1)}{s^{1/2-\eps}},\]
which is nonvanishing and $O(1)$ for $y=O(1)$.
For a fixed $s$ and with $a=U(s,\Gamma(s)+y)$ it follows that
\begin{eqnarray}
	\frac{d}{ds}U(s,\Gamma(s)+y) 
	&=& U_\eta(s,\Gamma(s)+y)\big[X_s(s;a)-\dot\Gamma(s)\big]\nonumber\\
	&=& O(1)\big[X_s(s;a)-\dot\Gamma(s)\big]  \label{dU/ds},
\end{eqnarray}
for $y=O(1)$.
On the other hand, using \eq{U_def} and Lemma \ref{dZ/ds}, we get
\[\frac{d}{ds}U(s,\Gamma(s)+y) = \frac{\frac{d}{ds}Z(s,\Gamma(s)+y+1)}{Z(s,\Gamma(s)+y+1)}
-\frac{\frac{d}{ds}Z(s,\Gamma(s)+y)}{Z(s,\Gamma(s)+y)}=\frac{O(1)}{s^{3/2-\eps}},\]
such that
\beq\label{X_dot-Gamma_dot}
	X_s(s;a)-\dot\Gamma(s)=\frac{O(1)}{s^{3/2-\eps}}.
\eeq
Transforming back to $t$ and $x$ coordinates, 
it follows that the $a$-level curve $x(t;a)$ of $u(t,x)$ satisfies
\beq\label{x_dot}
	x_t(t;a)=\dot\gamma(t)+\frac{O(1)}{(-t)^{3/2-\eps}},
\eeq
for $a$ sufficiently close to $a_0$.
We conclude that there exists a smooth function $c(a)$ defined on an interval about $a_0$ such that 
\beq\label{x}
	x(t;a)=c(a)+\gamma(t)+\frac{O(1)}{(-t)^{1/2-\eps}},
\eeq
for $t\in\,[\,-2T,0\,]$.

\section{ Estimates on the total variation}\label{total_variation}
We now fix the solution $u(t,x):=U(-t,-x)$ defined in 
\eq{U_def} and set
\[u_j^n:=u(n,j).\]
By Lemma \ref{Cole-Hopf} these values provide
a solution to the upwind scheme \eq{u_godunov} for the particular flux
function $f(u)=\ln(\nu+\mu e^u)$. Recalling that we consider negative times we
let $v_j^n$ denote the corresponding solution of \eq{v_godunov} obtained
by prescribing vanishing $v$-data at time $-2T$, and our goal is 
to show that the total variation of $V(j):=v^{1}_{j+1}$ over the set
$\{j\in\mathbb Z\,|\, -T\leq j \leq 0\}$ is of order $\ln T$. Here $T$ is as in Section 
\ref{special_discr_appr}.
This will be accomplished in a series of lemmas which reduce the problem 
to an estimate involving the heat kernel and the curve $\gamma(t)$. 
In the course of doing this we will specify the precise assumptions on 
the function $g(u)$ in \eq{v_godunov}. 

Recalling that we have fixed $\lambda=1/2$ in \eq{v_godunov} we get that 
the discrete Green kernel $K^n_{k}$ 
for (the linear part of) \eq{v_godunov} is given by 
\[K^n_k:=\frac{1}{2^n}\genfrac(){0cm}{0}nk.\]
With vanishing data at time $-2T$ we get that the solution at time 
step $n=1$ is given as
\beq\label{V}
	V(j):= v^{1}_{j+1}=\sum_{n=0}^{2T}\sum_{k\in\mathbb Z} \psi^{-n}_k K^{n}_{j-k-1},
\eeq
where we have set
\[\psi^{n}_k:=-\big[g(u_k^{n})-g(u_{k-1}^{n})\big]=-\int_{u(n,k-1)}^{u(n,k)}g'(\xi)\, d\xi.\]
We will make use of the following representation formula. 
\begin{lemma}\label{repr}
	Let $x(t;\xi)$ denote the $\xi$-level curve of the particular solution $u(t,x)$ specified above,
	and assume that the function $g(u)$ in \eq{v_godunov} has compact support.
	We then have the representation
	\beq\label{V_repr}
		V (j) =\int_0^\infty g'(\xi)\Big[\sum_{n=0}^{2T} K^{n}_{j-[\![x(-n;\xi)]\!]}\Big]\, d\xi.
	\eeq
\end{lemma}
\begin{proof}
	Using that the constructed solution $u(t,x)$ is monotone decreasing at each time we get that,
	\begin{eqnarray*}
		V(j) &=& \sum_{n=0}^{2T} \sum_{k\in\mathbb Z} \psi^{-n}_k K^{n}_{j-k+1}\\
		&=& -\sum_{n=0}^{2T} \sum_{k\in\mathbb Z} \left\{ \int_{u(-n,k-1)}^{u(-n,k)}g'(\xi)\, d\xi \right\} 
		K^{n}_{j-(k-1)}\\
		&=& -\sum_{n=0}^{2T} \sum_{k\in\mathbb Z} \int_{u(-n,k-1)}^{u(-n,k)}
		g'(\xi)K^{n}_{j-[\![x(-n;\xi)]\!]}\, d\xi \\
		&=& -\sum_{n=0}^{2T} \int_{u(-n,-\infty)}^{u(-n,+\infty)}g'(\xi)K^{n}_{j-[\![x(-n;\xi)]\!]}\, d\xi \\
		&=& \int_0^\infty g'(\xi)\Big[\sum_{n=0}^{2T} K^{n}_{j-[\![x(-n;\xi)]\!]} \Big]\, d\xi
	\end{eqnarray*}
\end{proof}
\noindent In the following computations we will use that the level curves $x(t;\xi)$ of $u(t,x)$ are 
well approximated by translates of $\gamma(t)$, which is true for large (negative) values 
of $t$ by \eq{x}. Now, for $T$ sufficiently large, the values $V(j)$ for $j\leq -\sqrt{T}$ are only
influenced exponentially little by values of $x(t;\xi)$ and $u(t,x)$ for  $t\in [\,-\sqrt{T},0\,]$. 
Furthermore, for $j \in [\,-T, -\sqrt{T}\, ]$, 
the $n$-summations in \eq{V} and \eq{V_repr} are unaffected, up to errors 
which are summable in $j$, by restricting $n$ to the interval 
\[I(j):= \Big[\frac{|j|}{\beta}-|j|^{1/2+\ve},\frac{|j|}{\beta}+|j|^{1/2+\ve}\Big],\]
where $0<\ve\ll 1$ and $\beta=\lambda_0-1/2>0$. 
This fact will be used repeatedly in the rest of this section.
From now on all values $j$ will be in $[\,-T, -\sqrt{T}\, ]$, i.e. far downstream.

For two functions $A(j)$ and $B(j)$ defined for integers $j$ we write $A \approx B$ to 
indicate that the difference of the functions  satisfy
\[\sum_{j\leq -\sqrt{T}}\left|A(j)-B(j)\right|<\infty.\] 
With $g$ as in Lemma \ref{repr} we thus have
\[V (j) \approx \int_0^\infty g'(\xi)\Big[\sum_{n\in I(j)} K^{n}_{j-[\![x(-n;\xi)]\!]}\Big]\, d\xi.\]
Since $n\in I(j)$ implies $|2(j-[\![x(-n;\xi)]\!])-n|\leq O(1) n^{1/2+\eps}$
we can use Lemma \ref{estimate} in the appendix  to conclude that  
\begin{eqnarray}
	&& V(j)-V(j-1) \approx \int_0^\infty\! g'(\xi)\Big[\sum_{n\in I(j)} K^n_{j-[\![x(-n;\xi)]\!]}-
	K^n_{j-[\![x(-n;\xi)]\!]-1}\Big]\, d\xi\nonumber \\\nonumber \\
	&&\quad\qquad\quad = 4\int_0^\infty g'(\xi)
	\Big[\sum_{n\in I(j)}G_x\big(\textstyle\frac{n}{2},2(j-[\![x(-n;\xi)]\!])-n\big)\Big]
	\, d\xi\nonumber \\
	&&\quad\qquad\qquad -4\int_0^\infty g'(\xi)
	\Big[\sum_{n\in I(j)}G_{xx}({\textstyle\frac{n}{2}},2(j-[\![x(-n;\xi)]\!])-n)\Big]\, d\xi\nonumber \\
	&&\quad\qquad\qquad + \int_0^\infty g'(\xi)
	\sum_{n\in I(j)}O(n^{-2+\ve})\,d\xi\nonumber \\\nonumber \\
	&&\quad\qquad\quad \approx 4\int_0^\infty g'(\xi)
	\Big[\sum_{n\in I(j)}G_x\big(\textstyle\frac{n}{2},2(j-[\![x(-n;\xi)]\!])-n\big)\Big]\, d\xi
	\nonumber\\
	&&\qquad\qquad\quad -4\int_0^\infty\! g'(\xi)
	\Big[\sum_{n\in I(j)}G_{xx}\big(\textstyle\frac{n}{2},2(j-[\![x(-n;\xi)]\!])-n\big)\Big]\, d\xi
	\nonumber \\\nonumber \\
	&&\qquad\qquad\quad =: \mathcal A(j) - \mathcal B(j), 
	\label{two_terms}
\end{eqnarray}
where we have used that $|I(j)|=O(1)|j|^{1/2+\ve}$ and that $n=O(|j|)$ for $n\in I(j)$.
We will estimate the variation of $V(j)$ on $[\,-T, -\sqrt{T}\, ]$ by considering 
$\mathcal A(j)$ and $\mathcal B(j)$ separately.

\subsection{Variation  of $\mathcal B(j)$} We start by defining the functions
\begin{eqnarray}
	H(y;\xi)&:= &\sum_{n\in I(j)}G_{xx}\big(\textstyle\frac{n}{2},2(y-[\![x(-n;\xi)]\!])-n\big),\label{big_H}\\
	h(y;\xi)&:= &\int_0^\infty G_{xx}\big(\textstyle\frac{t}{2},2(y-x(-t;\xi))-t\big)\, dt.\label{h}
\end{eqnarray}
\begin{lemma}\label{H=h}
	We have
	\beq
	H(j;\xi)\approx h(j;\xi)\quad\mbox{ uniformly for $\xi$ in compacts.}
	\eeq
\end{lemma}
\begin{proof}
	We recall the notation $h_1$ for the sawtooth function introduced in Section \ref{moving_source}. 
	Using the decay properties of the heat kernel, we have
	\begin{eqnarray*}
		&&\nqquad \big|H(j;\xi) - h(j;\xi)\big| = 
		\Big|\sum_{n\in I(j)}G_{xx}\big(\textstyle\frac{n}{2},2(j-[\![x(-n;\xi)]\!])-n\big)\\
		&&\qquad\qquad\qquad\qquad\qquad
		- \int_0^\infty G_{xx}\big(\textstyle\frac{t}{2},2(j-x(-t;\xi))-t\big)\, dt\Big|\\
		&=& \Big|\sum_{n\in I(j)}\Big[G_{xx}\big(\textstyle\frac{n}{2},2(j-x(-n;\xi))-n\big)
		+O(1) G_{xxx}\big(\textstyle\frac{n}{2},\tilde y(j,n;\xi)\big)\Big]\\
		&&\qquad\qquad\qquad\qquad
		- \int_0^\infty G_{xx}\big(\textstyle\frac{t}{2},2(j-x(-t;\xi))-t\big)\, dt\Big|\\		
		&\approx& \Big|\int_0^\infty G_{xx}\big(\textstyle\frac{t}{2},2(j-x(-t;\xi))-t\big)
		\big(1+h_1'(t)\big)\, dt\\
		&&\qquad\qquad\qquad\qquad\qquad 
		- \int_0^\infty G_{xx}\big(\textstyle\frac{t}{2},2(j-x(-t;\xi))-t\big)\, dt\Big|\\
		&=& \Big| \int_0^\infty \Big[{\textstyle\frac{1}{2}}G_{xxt}\big(\textstyle\frac{t}{2},2(j-x(-t;\xi))-t\big)\\
		&&\qquad \qquad \quad +2[x_t(-t;\xi)-{\textstyle\frac{1}{2}}]
		G_{xxx}\big(\textstyle\frac{t}{2},2(j-x(-t;\xi))-t\big)\Big]
		h_1(t)\, dt \Big|,
	\end{eqnarray*}
	where $\tilde y(j,n;\xi)$ is between  $2(j-x(n;\xi))-n$ and $2(j-[\![x(n;\xi)]\!])-n$.
	Using the decay properties of the heat  kernel, that $x_t=O(1)$,  and the fact that $t=O(1)|j|$ for $t\in I(j)$, 
	we conclude that 
	\[\big|H(j;\xi)- h(j;\xi)\big| \approx O(1)\int_{I(j)} \frac{1}{t^{5/2}}+\frac{1}{t^2}\, dt
	=\frac{O(1)}{|j|^{3/2-\ve}}.\]
\end{proof}
\noindent Note that the $O(1)$ in the last term here is uniform with respect to $\xi$ as $\xi$ varies over 
a compact set. Since $g$ is assumed to have compact support we conclude that the second term in 
\eq{two_terms} satisfies 
\[\int_0^\infty g'(\xi)\Big[\sum_{n\in I(j)}G_{xx}\big({\textstyle\frac{n}{2}},2(j-[\![x(n;\xi)]\!])-n\big)\Big]\, d\xi
\approx \int_0^\infty g'(\xi)h(j;\xi)\, d\xi.\]
Using \eq{x_dot} we now calculate
\begin{eqnarray*}
	h(j;\xi) &=& \int_0^\infty G_{xx}\big(\textstyle\frac{t}{2},2(j-x(-t;\xi))-t\big)\, dt\\
	&=& \int_0^\infty G_{xx}\big(\textstyle\frac{t}{2},2(j-x(-t;\xi))-t\big)[2x_t(-t;\xi)-1]\, dt \\
		&& + 2\int_0^\infty G_{xx}\big(\textstyle\frac{t}{2},2(j-x(-t;\xi))-t\big)[1-x_t(-t;\xi)]\, dt \\
		&=& \int_0^\infty G_{xx}\big(\textstyle\frac{t}{2},2(j-x(-t;\xi))-t\big)[2x_t(-t;\xi)-1] \\
		&&\qquad\qquad\qquad \qquad\qquad\qquad + 
		\frac{1}{2}G_{xt}\big(\textstyle\frac{t}{2},2(j-x(-t;\xi))-t\big)\, dt \\
		&& +2\int_0^\infty G_{xx}\big({\textstyle\frac{t}{2}},2(j-x(-t;\xi))-t\big)
		\Big((1-\lambda_0)-\frac{1}{\sqrt{t}}+\frac{O(1)}{t^{3/2}}\Big)\, dt \\
		&& -\frac{1}{2}\int_0^\infty G_{xxx}\big(\textstyle\frac{t}{2},2(j-x(-t;\xi))-t\big)\, dt,
\end{eqnarray*}
where we have added and subtracted $G_{xt}=G_{xxx}$. The first integral on the right 
vanishes so that rearranging gives,
\begin{eqnarray}
	(2\lambda_0-1) h(j;\xi)&\approx& -2\int_0^\infty G_{xx}
	\big({\textstyle\frac{t}{2}},2(j-x(-t;\xi))-t\big)\frac{dt}{\sqrt{t}} 
	+\int_{I(j)} \frac{O(1)}{|j|^2}\, dt\nonumber \\
	&=& \int_{I(j)} \frac{O(1)}{|j|^2}\, dt \, =\,  \frac{O(1)}{|j|^{3/2-\ve}}.
\end{eqnarray}
Hence $h(j;\xi)$ is summable in $j$, uniformly for $\xi$ in compacts.
It follows from Lemma \ref{H=h} that the same is true for $H(j;\xi)$, and we
conclude that the function 
$\mathcal B(j)$ is summable. Hence it contributes only a finite 
amount to the total variation of $V(j)$ independently of $T$.

\subsection{Variation of $\mathcal A(j)$} 
It remains to show that the variation of $\mathcal A(j)$ in \eq{two_terms} grows indefinitely 
with $T$. This requires a more detailed analysis and we start by 
defining the functions
\begin{eqnarray}
	L(y;\xi)&:= &\sum_{n\geq 0}G_{x}\big(\textstyle\frac{n}{2},2(y-x(-n;\xi))-n\big),\label{big_L}\\
	l(y;\xi)&:= &\int_0^\infty G_{x}\big(\textstyle\frac{t}{2},2(y-x(-t;\xi))-t\big)\, dt.\label{l}
\end{eqnarray}
Note that the arguments in neither of these functions involve integer parts.
\begin{lemma}\label{L=l}
	We have
	\beq
	L(j;\xi)\approx l(j;\xi)\quad\mbox{ uniformly for $\xi$ in compacts.}
	\eeq
\end{lemma}
\begin{proof}
	As above we use the function $h_1$ introduced in Section \ref{moving_source} and the 
	decay properties of the heat kernel:
	\begin{eqnarray*}
		&&\nqquad \big|L(j;\xi)-l(j;\xi)\big| 
		= \Big|\int_0^\infty G_{x}\big(\textstyle\frac{t}{2},2(j-x(-t;\xi))-t\big)
		\big(1+h_1'(t)\big)\, dt\\
		&&\qquad \qquad \qquad \qquad - \int_0^\infty 
		G_{x}\big(\textstyle\frac{t}{2},2(j-x(-t;\xi))-t\big)\, dt\Big|\\
		&=&\Big| \int_0^\infty \Big[{\textstyle\frac{1}{2}}
		G_{xt}\big(\textstyle\frac{t}{2},2(j-x(-t;\xi))-t\big)\\
		&&\qquad\qquad + [2x_t(-t;\xi)-1]
		G_{xx}\big(\textstyle\frac{t}{2},2(j-x(-t;\xi))-t\big)\Big]h_1(t)\, dt\Big|\\
		&\leq& \frac{O(1)}{|j|^{3/2-\ve}} +\Big|\int_0^\infty [2x_t(-t;\xi)-1]
		G_{xx}\big(\textstyle\frac{t}{2},2(j-x(-t;\xi))-t\big)h_1(t)\, dt\Big|.
	\end{eqnarray*}
	To estimate the last integral we apply \eq{x_dot} and the function $h_2$ introduced 
	in Section \ref{moving_source}, to get 
	\begin{eqnarray*}
		&&\int_0^\infty [2x_t(-t;\xi)-1]
		G_{xx}\big(\textstyle\frac{t}{2},2(j-x(-t;\xi))-t\big)h_1(t)\, dt \\ \\
		&&= \int_0^\infty \left(2\lambda_0-1 +\frac{O(1)}{\sqrt{t}}\right)
		G_{xx}\big(\textstyle\frac{t}{2},2(j-x(-t;\xi))-t\big)h_1(t)\, dt \\ \\
		&&\approx (2\lambda_0-1)  \int_0^\infty 
		G_{xx}\big({\textstyle\frac{t}{2}},2(j-x(-t;\xi))-t\big)h_2'(t)\, dt  +\frac{O(1)}{|j|^{3/2-\ve}}\\ \\
		&&=O(1)  \int_0^\infty G_{xxt}\big({\textstyle\frac{t}{2}},2(j-x(-t;\xi))-t\big)h_2(t)\, dt \\ \\
		&&\quad+O(1) \int_0^\infty G_{xxx}\big({\textstyle\frac{t}{2}},2(j-x(-t;\xi))-t\big)h_2(t)\, dt
		+\frac{O(1)}{|j|^{3/2-\ve}}\\ \\
		&&=\frac{O(1)}{|j|^{3/2-\ve}},
	\end{eqnarray*}
	and the conclusion follows.
\end{proof}
\begin{lemma}\label{l_int}
	$l(j;\xi)\approx 0$ uniformly for $\xi$ in compacts.
\end{lemma}
\begin{proof}
	Using \eq{x_dot} we have 
\begin{eqnarray*}
	l(j;\xi) &=&\int_0^\infty G_{x}\big(\textstyle\frac{t}{2},2(j-x(-t;\xi))-t\big)\, dt\\
	&=& \int_0^\infty G_{x}\big(\textstyle\frac{t}{2},2(j-x(-t;\xi))-t\big)[2x_t(-t;\xi)-1]\, dt \\
		&& + 2\int_0^\infty G_{x}\big(\textstyle\frac{t}{2},2(j-x(-t;\xi))-t\big)[1-x_t(-t;\xi)]\, dt \\
		&=& \int_0^\infty G_{x}\big(\textstyle\frac{t}{2},2(j-x(-t;\xi))-t\big)[2x_t(-t;\xi)-1] \\
		&&\qquad\qquad\qquad\qquad\qquad\qquad\quad + 
		{\textstyle\frac{1}{2}}G_{t}\big(\textstyle\frac{t}{2},2(j-x(-t;\xi))-t\big)\, dt \\ \\
		&&+2\int_0^\infty G_{x}\big({\textstyle\frac{t}{2}},2(j-x(-t;\xi))-t\big)
		\Big((1-\lambda_0)-\frac{1}{\sqrt{t}}+\frac{O(1)}{t^{3/2}}\Big)\, dt \\ \\
		&&-\frac{1}{2}h(j;\xi),
\end{eqnarray*}
where we have added and subtracted $G_{t}=G_{xx}$ and used \eq{h}. 
The first integral on the right-hand side here vanishes. Recalling the definition of 
$l(j;\xi)$, rearranging, and using the decay properties of the heat kernel together with 
what we have already proved about $h(j;\xi)$, we get
\[l(j;\xi) \approx \frac{2}{1-2\lambda_0}\int_0^\infty 
G_{x}\big({\textstyle\frac{t}{2}},2(j-x(-t;\xi))-t\big)\frac{dt}{\sqrt{t}}\, .\]
To estimate this last integral we repeat the same procedure:
\begin{eqnarray*}
	&&\nquad\int_0^\infty  G_{x}\big({\textstyle\frac{t}{2}},2(j-x(-t;\xi))-t\big)\frac{dt}{\sqrt{t}}\\ \\
	&&=\int_0^\infty  G_{x}\big({\textstyle\frac{t}{2}},2(j-x(-t;\xi))-t\big)\frac{[2x_t(-t;\xi)-1]}{\sqrt{t}}\, dt\\ \\
	&&\quad +2 \int_0^\infty  G_{x}\big({\textstyle\frac{t}{2}},2(j-x(-t;\xi))-t\big)\frac{[1-x_t(-t;\xi)]}{\sqrt{t}}\, dt\\ \\
	&&= \int_0^\infty  \Big\{[2x_t(-t;\xi)-1]G_{x}\big({\textstyle\frac{t}{2}},2(j-x(-t;\xi))-t\big)\\
	&&\qquad\qquad\qquad\qquad\qquad\qquad\quad  +
	{\textstyle\frac{1}{2}} G_{t}\big({\textstyle\frac{t}{2}},2(j-x(-t;\xi))-t\big)\Big\}\frac{dt}{\sqrt{t}}\\ \\
	&&\quad +2(1-\lambda_0) \int_0^\infty  G_{x}\big({\textstyle\frac{t}{2}},2(j-x(-t;\xi))-t\big)\frac{dt}{\sqrt{t}}\\ \\
	&&\quad +O(1)\int_0^\infty  G_{x}\big({\textstyle\frac{t}{2}},2(j-x(-t;\xi))-t\big)\frac{dt}{t}\\ \\
	&&\quad -\frac{1}{2} \int_0^\infty  G_{xx}\big({\textstyle\frac{t}{2}},2(j-x(-t;\xi))-t\big)\frac{dt}{\sqrt{t}}.
\end{eqnarray*}
Integration by parts and decay of $G(t)$ shows that the first integral here is $O(1)/|j|^{3/2-\ve}$.
Rearranging and using the decay properties of the heat kernel thus gives
\[\int_0^\infty  G_{x}\big({\textstyle\frac{t}{2}},2(j-x(-t;\xi))-t\big)\frac{dt}{\sqrt{t}}
=\frac{O(1)}{|j|^{3/2-\ve}}.\]
We conclude that 
\[l(j;\xi) \approx \frac{O(1)}{|j|^{3/2-\ve}},\]
which completes the proof.
\end{proof}
From Lemma \ref{L=l} and Lemma \ref{l_int} it follows that $L(j;\xi)$ is summable with respect to $j$,
uniformly for $\xi$ in compacts. We conclude that the function 
\beq\label{M}
	M(j;\xi) :=  \sum_{n\in I(j)}G_{x}\big({\textstyle\frac{n}{2}},2(j-[\![x(-n;\xi)]\!])-n\big)
\eeq
appearing in the first term in \eq{two_terms}, may be expanded as 
\begin{eqnarray*}
	M(j;\xi)&\approx& \sum_{n\geq 0} G_{x}\big({\textstyle\frac{n}{2}},2(j-[\![x(-n;\xi)]\!])-n\big)
	\, -\, L(j;\xi)\\
	&\approx& \sum_{n\geq 0} G_{xx}\big({\textstyle\frac{n}{2}},2(j-x(-n;\xi))-n\big)\big(\!\!\big(x(-n;\xi)\big)\!\!\big),
\end{eqnarray*}
where $\big(\!\!\big(a\big)\!\!\big):=a-[\![a]\!]$ denotes fractional part. Assume now that $\lambda_0$ is a rational 
number,
\[\lambda_0=\frac{p}{q},\]
say, and write $n=mq+i$ for $1\leq i\leq q$, to get that
\begin{eqnarray*}
	M(j;\xi) &\approx& \sum_{m\geq 0} \sum_{i=1}^q 
	G_{xx}\big({\textstyle\frac{mq+i}{2}},2(j-x(-mq-i;\xi))-mq-i\big)\\
	&&\qquad\qquad\qquad\qquad\qquad\qquad\quad \cdot \big(\!\!\big(x(-mq-i;\xi)\big)\!\!\big), 
\end{eqnarray*}
Taylor expanding $G_{xx}$ about the points 
\[\big({\textstyle\frac{mq}{2}},2(j-x(-mq;\xi))-mq\big),\] 
and using the decay properties of the heat kernel, yields 
\[M(j;\xi)\approx \sum_{m\geq 0} 
	G_{xx}\big({\textstyle\frac{mq}{2}},2(j-x(-mq;\xi))-mq\big)
	\left[ \sum_{i=1}^q\big(\!\!\big(x(-mq-i;\xi)\big)\!\!\big)\right].\]
Recalling from \eq{x} that the level curves $x(t;\xi)$ are given as 
\[x(-t;\xi)=c(\xi) -\lambda_0 t +2\sqrt{t}+\frac{O(1)}{t^{1/2-\eps}},\qquad \mbox{for $t>0$},\]
we obtain
\begin{eqnarray}
	M(j;\xi) &\approx&   \sum_{m\geq 0} 
	G_{xx}\big({\textstyle\frac{mq}{2}},2(j-x(-mq;\xi))-mq\big)
	\left[ \sum_{i=1}^q\big(\!\!\big(x(-mq-i;\xi)\big)\!\!\big)\right]\nonumber \\ \nonumber \\
	&=&  \sum_{m\geq 0} G_{xx}\big({\textstyle\frac{mq}{2}},2(j-x(-mq;\xi))-mq\big)\nonumber \\
	&&\qquad\qquad \cdot 
	\left[ \sum_{i=1}^q\big(\!\!\big(c(\xi) -\textstyle\frac{ip}{q} +2\sqrt{mq+i}+
	\frac{O(1)}{\sqrt{mq}}\big)\!\!\big)\right].\label{M1}
\end{eqnarray}
\begin{lemma}\label{M_approx}
	We have 
	\[M(j;\xi) \approx  \sum_{m\geq 0} G_{xx}\big({\textstyle\frac{mq}{2}},2(j-x(-mq;\xi))-mq \big)
	\left[\sum_{i=1}^q\big(\!\!\big(c(\xi) -\textstyle\frac{ip}{q} +2\sqrt{mq}\big)\!\!\big)\right]\]
	uniformly for $\xi$ in compacts.
\end{lemma}
\begin{proof}
	By \eq{M1} it suffices to show that
	\begin{eqnarray*}
	&&\nquad \sum_j\Big| \sum_{m\geq 0} 
	G_{xx}\big({\textstyle\frac{mq}{2}},2(j-x(-mq;\xi))-mq\big) \nonumber \\
	&&\quad \cdot \left[ \sum_{i=1}^q 
	\big(\!\!\big(c(\xi) -\textstyle\frac{ip}{q} +2\sqrt{mq}+\delta(m,i)\big)\!\!\big)-
	\big(\!\!\big(c(\xi) -\textstyle\frac{ip}{q} +2\sqrt{mq}\big)\!\!\big)\right]\Big|<\infty,
	\end{eqnarray*}
	where the summations are over those $m$ and $i$ for which $mq+i\in I(j)$.
	We have also set
	\[\delta(m,i):=2\sqrt{mq+i}-2\sqrt{mq}+\frac{O(1)}{\sqrt{mq}},\]
	which is of magnitude $O(1)/\sqrt{|j|}$ for all $m$ and $i$ in question. 
	Next, for fixed $\xi$ and $i$, let
	\[b:=c(\xi)-\textstyle\frac{ip}{q},\]
	and define 
	\[\alpha(s):=b+2\sqrt{ms},\qquad \beta(s):= b+2\sqrt{ms}+\delta(s,i).\]
	Now, as $s$ ranges over $I(j)$, both functions $\big(\!\!\big(\alpha(s)\big)\!\!\big)$ 
	and $\big(\!\!\big(\beta(s)\big)\!\!\big)$ jump at most an $O(1)|j|^\eps$ number of times.
	Also, if neither of these functions jump in an interval $(s-1,s+1)$, then they are 
	$O(1)/\sqrt{|j|}$ close on that interval. It follows that the expression above is dominated by 
	a sum of the form
	\[\sum_j\left\{ \frac{O(1)|j|^\eps}{|j|^{3/2}} + \frac{O(1)|j|^{1/2+\eps}}{|j|^{3/2}\sqrt{|j|}} \right\},\]
	which is finite.
\end{proof}
\noindent Making use of the identity
\[\sum_{i=1}^q\big(\!\!\big(a- \textstyle\frac{ip}{q} \big)\!\!\big)=\big(\!\!\big( qa\big)\!\!\big)+\frac{q-1}{2}\]
and what we have already shown for $H(j;\xi)$, we infer that
\begin{eqnarray*}
	M(j;\xi) &\approx&   \sum_{m\geq 0} 
	G_{xx}\big({\textstyle\frac{mq}{2}},2(j-x(-mq;\xi))-mq\big)\\
	&&\qquad\qquad\qquad\qquad\qquad \cdot
	\Big[\big(\!\!\big(qc(\xi) +2q\sqrt{mq}\big)\!\!\big)+\textstyle\frac{q-1}{2}\Big]\\\\
	&\approx& \sum_{m\geq 0} 
	G_{xx}\big({\textstyle\frac{mq}{2}},2(j-x(-mq;\xi))-mq\big)
	\big(\!\!\big(qc(\xi) +2q\sqrt{mq}\big)\!\!\big).
\end{eqnarray*}
We next want to approximate $M(j;\xi)$ with an integral and we 
proceed as above for $H(j;\xi)$ and $h(j;\xi)$. We define
\[N(j;\xi) :=  \int_0^\infty G_{xx}\big({\textstyle\frac{qs}{2}},2(j-x(-qs;\xi))-qs\big)
	\big(\!\!\big(qc(\xi) +2q\sqrt{qs}\big)\!\!\big)\, ds.\]
\begin{lemma}\label{M=N}
	 $M(j;\xi) \approx  N(j;\xi)$ uniformly for $\xi$ in compacts.
\end{lemma}
\begin{proof}
	Using the function $h_1$ and decay properties of the heat kernel we have 
	\begin{eqnarray*}
		&&\nquad M(j;\xi)-N(j;\xi) \\ \\
		&&\nquad \approx  \int_0^\infty G_{xx}\big({\textstyle\frac{qs}{2}},2(j-x(-qs;\xi))-qs\big)
		\big(\!\!\big(qc(\xi) +2q\sqrt{qs}\big)\!\!\big)h_1'(s)\, ds\\ \\
		&&\nquad= - \int_0^\infty \frac{d}{ds}\Big[G_{xx}\big({\textstyle\frac{qs}{2}},2(j-x(-qs;\xi))-qs\big)
		\big(\!\!\big(qc(\xi) +2q\sqrt{qs}\big)\!\!\big)\Big] h_1(s)\, ds\\ \\
		&&\nquad= - q\int_0^\infty \Big[ G_{xxx}\big({\textstyle\frac{qs}{2}},2(j-x(-qs;\xi))-qs\big)
		[2x_t(-qs;\xi)-1]\\ \\
		&&\qquad\qquad
		+ {\textstyle\frac{1}{2}}G_{xxt}\big({\textstyle\frac{qs}{2}},2(j-x(-qs;\xi))-qs\big) \Big]
		\big(\!\!\big(qc(\xi) +2q\sqrt{qs}\big)\!\!\big) h_1(s)\, ds\\ \\
		&& -\int_0^\infty G_{xx}\big({\textstyle\frac{qs}{2}},2(j-x(-qs;\xi))-qs\big)
		\left\{ \frac{d}{ds} \big(\!\!\big(qc(\xi) +2q\sqrt{qs}\big)\!\!\big)\right\}h_1(s) \, ds \\ \\
		&&\nquad \approx  -\int_0^\infty G_{xx}\big({\textstyle\frac{qs}{2}},2(j-x(-qs;\xi))-qs\big)
		\Big\{ {\textstyle\frac{q^{3/2}}{\sqrt{s}}} -\sum_{k\geq 0} \delta_{s_k}(s) \Big\}h_1(s) \, ds \\ \\
		&&\nquad \approx  \sum_{\{k\geq 0\, |\, s_k\in I(j)\}} 
		G_{xx}\big({\textstyle\frac{qs_k}{2}},2(j-x(-qs_k;\xi))-qs_k\big)h_1(s_k),
	\end{eqnarray*}
	where the $s_k$ denote the $s$-values for which $qc(\xi) +2q\sqrt{qs}$ is an integer.
	Clearly $s_k=O(k^2)$, such that the last sum contains only $O(1)|j|^{1/4+\eps}$ terms.
	Since each term in this sum is $O(|j|^{-3/2})$ the conclusion follows.
\end{proof}
\noindent 
Making a change of variables ($\tau=qs$) and restricting the integration to $I(j)$
we get that 
\beq\label{M_1}
	M(j;\xi)\approx \frac{1}{q} \int_{I(j)} G_{xx}\big({\textstyle\frac{\tau}{2}},2(j-x(-\tau;\xi))-\tau\big)
	\big(\!\!\big(qc(\xi) +2q\sqrt{\tau}\big)\!\!\big)\, d\tau.
\eeq
Recalling that the level curve $x(-\tau;\xi)$ is given as
\beq\label{c}
	x(-\tau;\xi)=c(\xi)-\lambda_0\tau+2\sqrt{\tau}+\frac{O(1)}{\sqrt{\tau}},
\eeq
we make a further change of variables $\tau \mapsto \si$ where
\[\tau=\frac{|j|}{\beta}+\frac{\sqrt{|j|}}{\beta^{3/2}}(\frac{\si}{\sqrt{2}}-2).\]
A straightforward calculation then yields
\[G_{xx}\big({\textstyle\frac{\tau}{2}},2(j-x(-\tau;\xi))-\tau\big)=
\frac{C}{|j|^{3/2}}\big(2\si^2-1\big)e^{-\si^2}+\frac{O(1)}{|j|^{2-\delta}},\]
which holds for $\tau\in I(j)$, or, equivalently, for $|\si-2\sqrt{2}|\leq O(1)|j|^\eps$. Here, and below, $C$ 
denotes various explicit numerical constants. 
Substituting into \eq{M_1} we conclude that $M(j;\xi)$ satisfies
\begin{eqnarray*}
	M(j;\xi) &\approx& \frac{C}{|j|}\int_{-\infty}^{\infty} (2\si^2-1)e^{-\si^2}\\
	&&\qquad\qquad\qquad \cdot
	\left(\!\!\!\left( qc(\xi) + 2q\sqrt{\textstyle\frac{|j|}{\beta}} 
	+ {\textstyle\frac{q(\si-\sqrt{8})}{\sqrt{2}\beta}}
	+{\textstyle\frac{O(1)}{|j|^{1/2-\eps}}}\right)\!\!\!\right)\, d\si \\ \\
	&\approx& \frac{C}{|j|}\int_{-\infty}^{\infty} (2\si^2-1)e^{-\si^2}
	\left(\!\!\!\left( qc(\xi) + 2q\sqrt{\textstyle\frac{|j|}{\beta}} 
	+ {\textstyle\frac{q(\si-\sqrt{8})}{\sqrt{2}\beta}}\right)\!\!\!\right)\, d\si.
\end{eqnarray*}

\medskip 

We can now finally proceed to estimate the variation of $V(j)$. 
We recall that the term $\mathcal B(j)$ in \eq{two_terms} is summable 
in $j$. Using the expression above for $M(j;\xi)$ we can estimate the variation 
of the term $\mathcal A(j)$ in \eq{two_terms}. 
Setting 
\beq\label{z}
	z_j:=2q\sqrt{\textstyle\frac{|j|}{\beta}}
\eeq 
we thus get that 
\begin{eqnarray*}
	&&V(j)-V(j-1) \approx 4\int_0^\infty g'(\xi)M(j;\xi)\, d\xi\\ \\
	&&\quad = \frac{C}{|j|}\int_{-\infty}^{\infty} (2\si^2-1)e^{-\si^2}
	\left\{\int_0^\infty g'(\xi)\left(\!\!\left( qc(\xi) + 
	z_j +{\textstyle\frac{q(\si-\sqrt{8})}{\sqrt{2}\beta}}\right)\!\!\right)\, d\xi \right\}\, d\si\\ \\
	&&\quad = \frac{C}{|j|}\int_{-\infty}^{\infty} (2\si^2-1)e^{-\si^2}\pi(\si;z_j)\, d\si\, .
\end{eqnarray*}
Here
\beq\label{5.16}\pi(\si;z):=\int_0^\infty g'(\xi)\left(\!\!\left( qc(\xi) + 
	z +{\textstyle\frac{q(\si-\sqrt{8})}
{\sqrt{2}\beta}}\right)\!\!\right)\, d\xi\, ,\eeq
with $c(\xi)$ as in \eq{c}. Defining the function 
\beq\label{5.17}
\Pi(z):=\int_{-\infty}^{\infty} (2\si^2-1)e^{-\si^2}\pi(\si;z)\, d\si\,,
\eeq
we observe that $\Pi(z)$ is periodic with period $1$ and 
(since $z_{-j}=z_j$) that
\beq\label{tv}
	\underset{-T\leq j\leq -\sqrt{T}}{\rm Tot. Var. }\,\, V(j)\,\,  \geq\,\, 
	\sum_{\sqrt{T}\leq j\leq T} \frac{C|\Pi(z_j)|}{j} \,\,-\,\,C_0,
\eeq
for some finite constant $C_0$ independent of $T$. 

In order to complete the argument showing that $V$ has large total variation 
we need to introduce the following condition. 

\medskip

\noindent
{\em {\bf Assumption (A):} The smooth function $g$ has compact support 
and is such that the 
$1$-periodic function $\Pi(z)$ at \eq{5.17} does not vanish identically.} 

\medskip
\noindent 
It is not difficult to see that 
this condition is satisfied for a large class of functions 
$g$. For example, if $g'$ is formally replaced by 
a Dirac delta-function $\delta_{u^*}$ 
concentrated at a point 
$u^*$ close to $a_0$, then we can directly compute
\[\Pi(z) = \int_{-\infty}^{\infty} (2\si^2-1)e^{-\si^2}
\left(\!\!\left( qc(u^*) + 
z +{\textstyle\frac{q(\si-\sqrt{8})}{\sqrt{2}\beta}}
\right)\!\!\right)\, d\sigma\, \not\equiv 0.\] 
By continuity, any function $g$ whose derivative approximates
$\delta_{u^*}$ will still satisfy the Assumption (A).

We now fix a smooth function $g$ which satisfies (A),
and observe that the resulting function $\Pi$ is continuous. Therefore, for 
a sufficiently small number $\eps$, $0<\eps<\sup_z |\Pi(z)|$, there exist numbers
$a< b$ in $(0,1)$ such that
\[|\Pi(z)| \geq \eps \qquad\forall z\in (n+a,n+b)\, \, =:J_n,\] 
for all $n\geq 0$.
Recalling the definition \eq{z} of $z_j$ we define, for each $n\in\mathbb N$, the integers 
\[j(n):=\quad\mbox{smallest integer $j$ for which $z_j\in J_n$\, },\]
and
\[k(n):=\quad\mbox{largest integer $j$ for which $z_j\in J_n$\, }.\]
That is,
\[j(n)=\left\lceil\frac{\beta (n+a)^2}{4q^2}\right\rceil \qquad\mbox{and}\qquad 
k(n)=\left\lfloor\frac{\beta(n+b)^2}{4q^2}\right\rfloor. \]
We conclude from \eq{tv} that
\begin{eqnarray*}
	\sum_{\sqrt{T}\leq j\leq T} \frac{|\Pi(z_j)|}{j} &\geq&  
	\sum_{\{n|\, J_n\subset [\sqrt{T},T]\}} 
\sum _{\{j\, |\, z_j\in J_n\}} \frac{|\Pi(z_j)|}{j} \\
	& \geq&  \sum_{\{n|\, J_n\subset [\sqrt{T},T]\}} 
\frac{\eps(k(n)-j(n))}{k(n)} \\
	&\geq&  \sum_{\{n|\, J_n\subset [\sqrt{T},T]\}} \frac{C}{n}
\quad \geq \quad C\ln T.
\end{eqnarray*}
This completes the proof and shows that there is no {\em a priori} 
bound available on the 
total variation of an approximate solution computed using the upwind scheme.

\section{Concluding remarks}\label{remarks}
The analysis shows that it is possible to
prescribe data for a strictly hyperbolic system of conservation laws 
in such a way that the total variation of the solution generated 
by the Godunov scheme increases by an arbitrarily large amount,
after a large  number of time steps. 
We observe that the total change in  the speed of the 
curve $\gamma$ decreases to zero as $T\to \infty$.
Therefore, as initial data for the $u$-component, we can take
an arbitrarily small perturbation of a single shock.

We summarize here the main features
of our example:
\v
\noindent{\bf 1.} Our $2\times 2$ system is strictly hyperbolic. 
One characteristic field is linearly degenerate, the other is
genuinely nonlinear.  The characteristic speeds are both 
strictly contained inside the interval $[0,1]$, so that the
usual linearized stability conditions hold.
\v
\noindent{\bf 2.} We can choose a sequence of initial data of the form 
$$u_\nu(0,x)=\left\{
\begin{array}{ll}
u^-+\phi_\nu(x)\qquad &\hbox{if}\quad x<0\,,\\
~~~u^+\qquad &\hbox{if}\quad x>0\,,\\
\end{array}
\qquad\qquad
v_\nu(0,x)=0\,,\right.$$
where the functions $\phi_\nu$ are smooth and satisfy
$$
\mbox{Tot. Var. } \{ \phi_\nu\}\to 0\,,
\qquad\qquad \|\phi_\nu\|_{{\mathcal C}^k}\to 0$$ 
for all $k$, as $\nu\to\infty$.   
For these initial data, the exact solutions have uniformly 
bounded total variation, but the corresponding
Godunov approximations $(u_\nu,\,v_\nu)$
satisfy
$$\mbox{Tot. Var. } \big\{v_\nu(T_\nu,\cdot)\big\}
\to\infty
\,,$$
for some sequence of times $T_\nu\to\infty$.
\v
\noindent{\bf 3.} Consider a second sequence of Godunov solutions
$$(\tilde u_\nu\,,~\tilde v_\nu)(t,x)= (u_\nu\,,~v_\nu)
(t\,,~x-1)\,.$$
By the previous estimates on the total variation,
it trivially follows that
$$\lim_{\nu\to\infty} \frac{\big\|v_\nu(T_\nu)-\tilde v_\nu(T_\nu)\big\|_{
{\bf L}^1}}
{\big\|u_\nu(0)-\tilde u_\nu(0)\big\|_{
{\bf L}^1}} = \infty\,.$$
In other words, the Godunov approximations are also 
unstable in the ${\bf L}^1$
norm, w.r.t.~perturbations of the initial data. This is in
sharp contrast with the stability of Glimm or front tracking
approximations \cite {LiuYang99}, \cite{BressanLiuYang99}.
\v
Our counterexample has many special features: 
it can be solved in triangular form, the first equation
contains
a very particular flux function, 
and the second equation is linear.
All these additional features allow us to 
perform explicit calculations
and derive rigorous estimates, but none of them
seems to be essential for the validity of the result.
Indeed, we expect that the instability highlighted 
in this paper will be a common feature of all
discrete schemes for systems of conservation laws, 
under generic conditions on the flux functions.
In this direction, we observe that 
in our example the oscillations in the $v$-component
are spread out on a very large spatial interval.
Requiring that the first field be genuinely nonlinear 
(instead of linearly degenerate) would not achieve
any appreciable decay in these oscillations.

We observe also that our results do {\em not} show that the upwind scheme 
fails to converge. In light of the convergence results in \cite{DingChenLuo89} 
(obtained through compensated compactness arguments) we expect that
one can in fact establish strong convergence also for the systems considered 
above. In this connection it is worth noticing that triangular $2\times 2$ 
systems of the form \eq{u_eq}-\eq{v_eq} are endowed 
with strictly convex entropies.

A somewhat subtle point in the preceding analysis is the issue of small data. 
As presented above, our examples do not immediately fit the 
setting of Glimm's theorem.
We started by giving a particular shock-wave solution of the 
first equation \eq{u_eq}, 
and {\em then} chose a flux $g(u)$ for the second equation. 
However it is now possible to 
go back and see that if we first fix a smooth function $g(u)$, then
we can give shock-wave solutions of \eq{u_eq} of arbitrarily small 
amplitude and 
with the properties: (i) its level curves travels with speed 
$\approx [\mbox{rational}] + 
1/\sqrt{-t}$, and (ii) the assumption (A) is satisfied. 
One can then redo the
previous 
analysis for the resulting system/data pair. Since the data for 
the second equation are 
identically zero, this provides an example within the setting of 
Glimm's theorem
for small BV data.

\section{Appendix}
\paragraph{Approximation of binomial coefficients}
For completeness we include the argument for the approximation 
of the discrete Green's function in terms of the heat kernel. Following 
the notation in Feller's book \cite{Feller68},
\[a_k(\nu):=\left(\frac{1}{2}\right)^{2\nu}\genfrac(){0cm}{0}{2\nu}{\nu+k},\]
we have
\begin{eqnarray*}
	a_k(\nu)&=&
	=a_0\frac{\exp \left[ \sum_{j=-k+1}^0 \ln \left(1-\textstyle\frac{j}{\nu}\right) \right] }
	{\exp \left[ \sum_{j=1}^k \ln \left(1-\textstyle\frac{j}{\nu}\right) \right]}\\\\
	&=& a_0\exp \left[-  \textstyle\frac{2}{\nu} \left(\sum_{j=1}^{k-1}j\right)
	 -\textstyle\frac{2}{3\nu^2} \left(\sum_{j=1}^{k-1}j^3\right)\right.\\ \\
	 && \qquad\qquad\left. - \textstyle\frac{k}{\nu} +  
	 \textstyle\frac{1}{2}\left(\textstyle\frac{k}{\nu}\right)^2	
	 -\textstyle\frac{1}{3}\left(\textstyle\frac{k}{\nu}\right)^3 
	 +\textstyle\frac{1}{4}\left(\textstyle\frac{k}{\nu}\right)^4
	 +O\left(\textstyle\frac{k^6}{\nu^5}\right)\right].
\end{eqnarray*}
Thus,
\[a_k(\nu)=a_0e^{-\frac{k^2}{\nu}}\exp\left[\textstyle\frac{1}{2}\left(\textstyle\frac{k}{\nu}\right)^2
-\textstyle\frac{1}{3}\left(\textstyle\frac{k}{\nu}\right)^3-
\textstyle\frac{(k-1)^2k^2}{6\nu^3}+\textstyle\frac{1}{4}\left(\textstyle\frac{k}{\nu}\right)^4
+O\left(\textstyle\frac{k^6}{\nu^5}\right)\right].\]
Defining
\[h:=\sqrt{\frac{2}{\nu}}\, ,\qquad\qquad \mathfrak{N}(x):=\frac{1}{\sqrt{2\pi}}e^{-x^2/2},\]
and recalling Stirling's formula
\[n!=\sqrt{2\pi}n^{n+1/2}e^{-n}\left(1+\frac{1}{12n}+\frac{1}{288n^2}
-\frac{139}{51840n^3}+\frac{O(1)}{n^4}\right),\]
we get
\[a_k(\nu)=h\mathfrak N(kh)e^{\eps_1-\eps_2}.\]
Here $\eps_1$ is due to the expansion above and $\eps_2$ is due to the 
approximation in Stirling's formula. We have
\[\eps_1=\eps_1(k,\nu):=\frac{1}{2}\left(\frac{k}{\nu}\right)^2
-\frac{1}{3}\left(\frac{k}{\nu}\right)^3-
\frac{(k-1)^2k^2}{6\nu^3}+O(1)\left(\frac{k}{\nu}\right)^4,\]
while 
\[\eps_2=\eps_2(\nu)=-\frac{1}{8\nu}+\frac{1}{192\nu^3}+\frac{O(1)}{\nu^4}.\]

\paragraph{Discrete Green's function}
Applying the approximation above to the discrete Green's function 
\[K^n_{k}:=\left(\frac{1}{2}\right)^{n}\genfrac(){0cm}{0}{n}{k},\]
we get for $n=2m,\, k=2l$ that
\begin{eqnarray*}
	K^n_{k}=a_{k-m}(m) &=& \sqrt{\frac{2}{m}}\mathfrak N \left( \sqrt{\textstyle\frac{2}{m}}(k-m) \right)e^{\eps_1-\eps_2}\\\\
	&=& 2G\left(\frac{n}{2},2k-n\right)e^{\eps_1-\eps_2},
\end{eqnarray*}
where $G(t,x)=e^{-x^2/4t}/2\sqrt{\pi t}$ denotes the heat kernel and 
\begin{eqnarray*}
	\eps_1&=&\eps_1(k-{\textstyle\frac{n}{2}}, {\textstyle\frac{n}{2}})
	= \frac{1}{2}\left(\frac{2k-n}{n}\right)^2-\frac{1}{3}\left(\frac{2k-n}{n}\right)^3\\\\
	&&\qquad\qquad\qquad-\frac{(2k-n-2)^2(2k-n)^2}{12n^3}+O\left(\left(\frac{2k-n}{n}\right)^4\right),
\end{eqnarray*}
and
\[\eps_2=\eps_2({\textstyle\frac{n}{2}})=\frac{1}{4n}-\frac{1}{24n^3}+\frac{O(1)}{n^4}.\]

\begin{lemma}\label{estimate}
	If $|n-2k|\leq O\left(n^{1/2+\delta} \right)\, (0<\delta\ll 1)$, then
	\[K^n_{k}-K^n_{k-1}=4\Big[G_x\left(\textstyle\frac{n}{2},2k-n\right)-
	G_{xx}\left(\textstyle\frac{n}{2},2k-n\right)\Big]+O\left(n^{-2+4\delta}\right).\]
\end{lemma}
\begin{proof}
	Letting
	\[\eps_1=\eps_1(k-{\textstyle\frac{n}{2}}, {\textstyle\frac{n}{2}}),
	\qquad \eps_2=\eps_2({\textstyle\frac{n}{2}}), \qquad \tilde\eps_1=\eps_1(k-n/2-1,n/2),\]
	we get that
	\begin{eqnarray*}
	K^n_{k}-K^n_{k-1} &=& 2\Big[G\left(\textstyle\frac{n}{2},2k-n\right)
	-G\left(\textstyle\frac{n}{2},2k-n-2\right)\Big]e^{\eps_1-\eps_2} \\\\
	&&+2G\left(\textstyle\frac{n}{2},2k-n-2\right)\left(e^{\eps_1-\eps_2}-e^{\tilde\eps_1-\eps_2}\right)\\\\
	&=& 4\left[G_x\left(\textstyle\frac{n}{2},2k-n\right)-
	G_{xx}\left(\textstyle\frac{n}{2},2k-n\right)+\frac{O(1)}{n^2}\right]e^{\eps_1-\eps_2} \\\\
	&&+2G\left(\textstyle\frac{n}{2},2k-n-2\right)\left(e^{\eps_1}-e^{\tilde\eps_1}\right)e^{-\eps_2}.
	\end{eqnarray*}
	Since $(2k-n)/n=O\left(n^{-1/2+\delta}\right)$ a calculation shows that 
	\[\eps_1-\eps_2=O\left(n^{-1+4\delta}\right),\]
	and that
	\[\eps_1-\tilde\eps_1=O\left(n^{-3/2+3\delta}\right).\]
	Thus
	\begin{eqnarray*}
	&&\nquad K^n_{k}-K_{n,k-1}\\\\
	&&= 4\left[G_x\left(\textstyle\frac{n}{2},2k-n\right)-
	G_{xx}\left(\textstyle\frac{n}{2},2k-n\right)+\frac{O(1)}{n^2}\right]\left(1+O\left(n^{-1+4\delta}\right)\right) \\\\
	&&\quad+O(n^{-1/2})\left(O\left(n^{-3/2+3\delta}\right)+O\left(n^{-5/2}\right)\right)\\\\
	&&=4\Big[G_x\left(\textstyle\frac{n}{2},2k-n\right)-
	G_{xx}\left(\textstyle\frac{n}{2},2k-n\right)\Big]+O(1)n^{-2+4\delta}.
	\end{eqnarray*}
\end{proof}

\paragraph{Acknowledgments} 
A. Bressan and P. Baiti:
Research supported in part by the European network HYKE, funded by the
EC as contract HPRN-CT-2002-00282, and by the italian MIUR-COFIN-2002
research project 2002017219 ``Equazioni Iperboliche e
Paraboliche Nonlineari".
H. K. Jenssen: Research supported in part by the NSF under grant DMS-0206631.


\end{document}